\documentclass[11pt]{article}

\makeatletter
\def\@citex[#1]#2{\if@filesw\immediate\write\@auxout
        {\string\citation{#2}}\fi
\def\@citea{}\@cite{\@for\@citeb:=#2\do
        {\@citea\def\@citea{,}\@ifundefined
        {b@\@citeb}{{\bf ?}\@warning
        {Citation `\@citeb' on page \thepage \space undefined}}
        {\csname b@\@citeb\endcsname}}}{#1}}
\newif\if@cghi
\def\cite{\@cghitrue\@ifnextchar [{\@tempswatrue
        \@citex}{\@tempswafalse\@citex[]}}
\def\citelow{\@cghifalse\@ifnextchar [{\@tempswatrue
        \@citex}{\@tempswafalse\@citex[]}}
\def\@cite#1#2{{\if@cghi\unskip$\null^{#1}$\else #1\fi\if@tempswa\typeout
        {warning: optional citation argument ignored: `#2'} \fi}}

\def\@biblabel#1{$\null^{#1}$}

\makeatother

\makeatletter

\@addtoreset{equation}{section}
\makeatother

\usepackage{graphics}  
\usepackage{amsmath,amssymb}  
\usepackage{epsf,amsfonts}
\usepackage{fontenc,indentfirst, delarray,epstopdf}
\DeclareGraphicsExtensions{eps}
\newtheorem{lemma}{Lemme}[section]
\newtheorem{proposition}[lemma]{Proposition}
\newtheorem{corollary}[lemma]{Corollary}
\newtheorem{definition}[lemma]{Definition}

\newcommand{\norm}[1]{\left\lVert#1\right\rVert}
\newcommand{\abs}[1]{\lvert#1\rvert}

\newcommand{\scl}[2]{\langle#1,#2\rangle}
\newcommand{\suup}[1]{ \underset{#1}{\sup} }

\newcommand{\tr}[1]{\text{Tr}\; #1}

\def\lp{\left(} 			
\def\rp{\right)}			
\def\dm{\lp\begin{array}} 		
\def\fm{\end{array}\rp}
\def\m2{M_2 \lp \cc \rp}	
\def\mn{M_n \lp \cc \rp}	
\def\re{\text{Re}}			

					\def\cc{{\mathbb{C}}}			\def\rr{{\mathbb{R}}}
\def\nn{{\mathbb{N}}}	\def\mm{{M}}			\def\zz{{\mathbb{Z}}} 		\def\ii{{\mathbb{I}}}
\def\aa{{\mathcal A}}						
\def\dd{{\mathcal D}}	\def\hh{{\mathcal H}}

\def\pp{{\mathcal P}}

\def\cinf{C^{\infty}\lp\mm\rp}

\def\ot{\otimes}
\def\xox{{\xi}_x}
\def\yox{{\xi}_y}
\def\xoz{{\zeta}_x}
\def\yoz{{\zeta}_y}

\begin{document} 


\title{\bf{Spectral distance on the circle}}     


\author{Pierre Martinetti
\\Dipartimento di fisica-GTC, \\Universit\`a di Roma "La
Sapienza", Piazzale Aldo Moro 5, I-00185
Roma
}
\maketitle
\begin{abstract}
\footnote{Paper published in Journal of Functional Analysis. Original version is available 
at http://dx.doi.org/10.1016/j.jfa.2008.07.018} A building block of noncommutative geometry is the observation that
most of the geometric information of a compact riemannian spin
manifold $M$ is encoded within its Dirac operator $D$. Especially
via Connes' distance formula one is able to extract from the spectral
properties of $D$ the  
geodesic distance on $M$. In this paper we investigate the distance $d$ encoded within
a covariant Dirac operator on a trivial $U(n)$-fiber bundle over the circle
with arbitrary connection.  It turns out that the
connected components of $d$ are tori whose dimension is given by the
holonomy of the connection. For $n= 2$ we explicitly compute $d$ on
all the connected components.  For $n\geq 2$ we restrict to a given fiber
and find that the distance is given by the trace of the module of a matrix. The latest is defined
by the holonomy and the coordinate of the points under consideration. This paper extends to arbitrary $n$ and arbitrary connection 
the results obtained in a previous work for $U(2)$-bundle with 
constant connection. It confirms interesting pro\-perties of the spectral distance with respect to another distance naturally associated to connection, namely the horizontal or Carnot-Carath\'eodory distance $d_H$.  Especially in case the connection has  irrational components, the connected components for $d$ are the closure of the connected components of $d_H$ within the euclidean topology on the torus. 
\end{abstract}
%

\section{Introduction}
\label{intro}  
In Connes framework a noncommutative geometry
consists in a {\it spectral triple}\cite{connes} $\aa ,\; \hh,\; D$
($\aa$ is an involutive algebra commutative or not, $\hh$ a Hilbert
space carrying a representation $\Pi$ of $\aa$ and $D$ a selfadjoint
operator on $\hh$ called Dirac operator) together with a graduation $\Gamma$ (called {\it
chirality} in physicists' words) and a real structure $J$ both
acting on $\hh$. In analogy with the commutative case  points are
recovered as pure states $\pp(\aa)$ of $\aa$ and the distance $d$
between two states $\omega$, $\omega'$ is defined as\cite{metrique}
\begin{equation}
\label{distance} d(\omega, \omega') \doteq \suup{a\in\aa}\left\{
\abs{\omega(a) - \omega'(a)}\,;\; \norm{[D,\Pi(a)]}\leq 1\right\}
\end{equation}
where the norm is the operator norm on $\hh$. 
This formula appears as a natural extension to the
noncommutative realm of the classical definition of the distance
 as
the length of the shortest path. Indeed  
to any compact, oriented, 
riemannian, spin manifold $\mm$ is naturally associated a spectral triple
\begin{equation}
\label{continu} \aa_E = \cinf,\; \hh_E=L_2(M,S),\; D_E=
-i\gamma^\mu\partial_\mu
\end{equation}
where $\hh_E$ is the space of square integrable spinors over $\mm$,
with functions acting by multiplication,   
and $D_E$ is the ordinary Dirac operator
of quantum field theory ($\gamma^\mu$ is the representation of the
Clifford algebra on $\hh_E$ and we omit the spin connection for it
commutes with $\Pi$) and it is not difficult to see that
(\ref{distance}) then coincides with the geodesic
distance defined by the riemannian structure of $\mm$. To get convinced
one
first identifies the pure states of $\cinf$ to the points of $M$ by Gelfand
duality, 
\begin{equation}
\label{gelfand}
\omega_x(f) = f(x),\quad  x\in \mm, f\in\cinf, \omega_x\in\pp(\cinf),
\end{equation}
then observes that $\norm{[D_E,\Pi(f)]}\leq 1$ makes 
$f$ be Lipschitz with constant $1$, so that the supremum of $\abs{\omega_x(f)-\omega_y(f)}$ in
(\ref{distance}) is
less than the geodesic distance between $x$ and $y$. The equality is
attained by considering
the function
\begin{equation}
\label{functionl}
L(z)\doteq d_{\text{geo}} (x,z)\quad \forall z\in\mm .
\end{equation}
Therefore
(\ref{distance}) is a reformulation of the classical definition of the
geodesic distance that 

- still makes sense in the noncommutative case;

- does not rely, in the commutative case, on any notion ill-defined
in a quantum framework such as trajectory
between points.

\noindent In fact only spectral
properties of $\aa$ and $D$ are involved, that
is why in the following we call $d$ the {\it spectral distance} associated to the
 geometry 
$(\aa, \hh, D)$. In previous papers we used the
terminology "noncommutative distance" but that might be confusing,
suggesting $d$ was not symmetric. Of course $d$ is symmetric and is
a distance in the strict mathematical sense.

Implementing a connection within a geometry $(\aa, \hh, D)$ amounts to
substituting $D$ with a {\emph covariant Dirac operator} $D_A$ (see
eq.(\ref{diraccov})). The
spectral distance becomes\footnote{Most of the time we omit the
symbol $\Pi$ and it should be clear from the context whether $a$
means an element of $\aa$ or its representation on $\hh$.} 
\begin{equation}
\label{distancefluc} d(\omega, \omega') \doteq \suup{a\in\aa}\left\{
\abs{\omega(a) - \omega'(a)}\,;\; \norm{[\mathcal{D}, a]}\leq 1\right\}
\end{equation}
where
$\mathcal{D}$ denotes the part of $D_A$ that does not commute with the
representation. Since there is no reason  for
$\left[\mathcal{D},a\right]$ to equal $\left[D,a\right]$,
(\ref{distancefluc}) differs from (\ref{distance}) and we talk of a {\emph{fluctuation of the metric}}.
 In this paper we study such fluctuations for
almost commutative geometries, that is to say products of the
continuous geometry (\ref{continu})
by another, yet undetermined, {\emph internal} geometry $(\aa_I,\; \hh_I,\;
D_I)$. The relevant algebra is then 
\begin{equation}
\aa =\aa_E \ot \aa_I = C^{\infty}(M,\aa_I)
\end{equation} represented on $\hh_E \ot \hh_I$ with a Dirac operator
$D$ given in (\ref{tripletprod}). In this case fluctuations of the metric amounts
to adding two pieces to $D$: 
a scalar field over $\mm$ denoted $H$ and a $1$-form field $\gamma^\mu
A_\mu$ (see eqs.(\ref{gaugg},\ref{gaugh})).
 When the internal algebra $\aa_I$ has finite dimension,
$A_\mu$ has value in the Lie algebra of the unitaries of $\aa_I$ and
provides the trivial $U(\aa_I)$-bundle of pure states,
\begin{equation}
\pp(\aa) = P\overset{\pi}{\rightarrow} \mm,
\end{equation}
with a connection. For this
reason $A_\mu$ is called the {\it gauge part} of the fluctuation.

 In [\citelow{kk}]
we have computed the fluctuation of the spectral distance
(\ref{distancefluc}) in almost commutative geometry for
scalar fluctuations only ($A_\mu = 0$) with various choices of
finite dimensional $\aa_I$. In [\citelow{cc}] we have studied pure
gauge fluctuations ($H=0$) for $\aa_I=\mn$ ($n\times n$ complex matrices). Then
the set of pure states $\pp(\aa)$ is a trivial $U(n)$-bundle $P$ on
$M$ with fiber $\pp(M_n(\cc)) \simeq \cc P^{n-1}$ equipped with a connection given locally by the $1$-form field $A_\mu$. 
On this bundle it turns out that the spectral distance (\ref{distancefluc}) is always
smaller or equal to another distance classically associated to a 
connection, namely the Carnot-Carath\'eodory (or horizontal or sub-riemannian)
distance\cite{montgomery} $d_H$. More about the horizontal distance is said in
section II, for the moment it is enough to underline that the
horizontal distance plays for gauge fluctuations the same
role as the geodesic distance in the commutative case, namely $d_H$
provides an upper bound in formula (\ref{distancefluc}) as
$d_{\text{geo}}$ provides an upper-bound in (\ref{distance}). In fact it was expected\cite{gravity} that the spectral
distance  and the horizontal distance were equal. But we pointed out in [\citelow{cc}] an
obstruction due to the holonomy of the connection that forbids $d$ to
equal $d_H$. More precisely when the
holonomy is trivial (i.e. the connection is flat
and $M$ is simply connected) then $d=d_H$ on all $P$. But when the holonomy is
non-trivial there no
longer exists an equivalent in $\aa$ of the function $L$ introduced in eq.(\ref{functionl}). In
other terms $d_H$ is an upper bound in eq.(\ref{distancefluc}) but not the
lowest upper bound. The obstruction is due to the self-intersecting points
of the minimal horizontal paths between $p, q\in P$, 
that is to say the points $p_i$  that project onto the
same point $z$ in the base  (see figure \ref{helixfig}). 
When there are too many such points one can prove that no equivalent
of the function $L$ exists. Therefore in order to compute the gauge
fluctuation of the metric in products of geometries one should first answer the
following question:\newline {\emph{do all the minimal
 horizontal curves between two points on a
fiber bundle have the same numbers of self-intersecting points ?}}\newline
Or,  more specifically:\newline
 {\emph{what is the minimum number
     of self-intersecting points
of a minimal
     horizontal path between two given points in a fiber bundle ?}}

 None of these question seems to have an answer in sub-riemannian
geometry\cite{montgomery}. If the base of the bundle has
sufficiently high-dimension, one may expect to deform a
minimal horizontal curve keeping its length fixed and reducing the
number of self-intersecting points. A way to escape the problem is to consider
a bundle with a low-dimensional base, for instance $S^1$. Then all horizontal
curves are 
horizontal-lifts of the circle
and for any two points in the bundle there exists at most one minimal
horizontal curve. Based on this example we found in [\citelow{cc}] that not only $d$ is
different from $d_H$, but that the two distances have properties rather
different: while $d_H$ is finite only between points that belong to
the same orbit of the holonomy group, $d$ has a larger connected
component. These differences have been
illustrated on a simple ($A_\mu = \text{constant}$) low dimensional
($\aa_I = M_2(\cc)$) example, focusing on
equatorial states (see definition, eq.(\ref{equatorial})
below). The present paper aims at generalizing this result to
arbitrary $A_\mu$ (non constant and of dimension $n\geq 2$), taking into account non-equatorial states. Although
self-contained, this work is better understood in relation with
[\citelow{cc}].
\newline

The paper organizes as follows: in section 2 we recall the basics
on products of spectral triples and fluctuations of the metric, as well
as the definition of the horizontal distance.  We briefly explain
why, due to the holonomy constraint, we restrict to  $U(n)$-bundle 
$P\overset{\pi}{\rightarrow} S^1$ over the circle and we introduce 
the main notations regarding the pure states $\left\{\xox; \xi\in\cc
    P^{n-1},x\in S^1\right\}$ of $C^\infty(S^1, \mn)$. Section 3 deals with
topology. We show that:

- the
connected components $\text{Con}(\xox)$ of the spectral distance on
$P$ with arbitrary
connection are tori $\mathbb{U}_\xi$ with dimension $n_c\leq n$
given by the holonomy; 

- the connected component $\text{Acc}(\xox)$
for the horizontal distance is a dense or a discrete subset of
$\mathbb{U}_\xi$, depending on the irrationality of the connection.

\noindent The main result is stated in proposition
\ref{propconneccomp} which extends to arbitrary connections and
dimensions the result of [\citelow{cc}]. In section 4 we explicitly compute the spectral
distance between any two states in $P$ for the case $n=2$, extending the results of [\citelow{cc}] to
non-equatorial states. Especially it turns out that the dependence
of the distance on the "altitude" $z_\xi$ of the states is far from
trivial, as explained in proposition \ref{propcarat} together with
its two related corollaries \ref{tplusmoinscorrol} and
\ref{corollcarat}. These rather technical results are discussed in
section 4.3. Section 5 deals with the $n>2$ case for which we
explicitly compute the spectral distance between states on the same
fiber. Quite nicely we find a simple expression in terms of the trace
of the module of a matrix whose elements are functions of the holonomy
and of the states under consideration.

\section{Why the circle ?}

This section does not contain new material. It aims at recalling the
basic features of the fluctuation of the metric in noncommutative
geometry, as well as the horizontal distance in a fiber bundle. We
then explain why, in order to study the relation between the
fluctuated spectral distance and the horizontal distance,  we restrict to bundles on the circle.

\subsection{Fluctuation of the metric in a product of geometries}

Connections are implemented within a geometry $(\aa, \hh,
\mathcal{D})$ by substituting $D$ with a so called {\it covariant
operator}
\begin{equation}
\label{diraccov} D_A \doteq D + A + JAJ^{-1}
\end{equation} where $A$ is a selfadjoint element of the set
$\Omega^1$ of 1-forms \footnote{We use Einstein summation over repeated indices in alternate
positions (up/down).}
\begin{equation}
\Omega^1 \doteq  \left\{a^i[D,b_i]\;;\; a^i,b_i\in\aa\right\}.
\end{equation}
For details on this statement the reader is referred to classical textbooks,
such as [\citelow{kob}] concerning the general theory of
connection and [\citelow{gravity}] for its implementation  within
noncommutative geometry. 
To our purpose it is sufficient
to know
that the properties defining
a spectral triple guarantees that $JAJ^{-1}$ commutes with the
representation. Thus only the part of $D_A$ that does not obviously commute with the
representation, namely
\begin{equation}
\label{tildeD} \mathcal{D} \doteq D + A,
\end{equation}
enters in the distance formula.

From now on we assume that $(\aa, \hh,D)$ is the product of the 
external spectral triple (\ref{continu}) by a so called internal
geometry $(\aa_I ,
\hh_I , D_I)$. Such a product, defined as
\begin{equation}
\label{tripletprod}
 \aa \doteq \aa_E\ot\aa_I,\; \hh\doteq \hh_E\ot
\hh_I,\; D\doteq D_E \ot \ii_I + \gamma^5\ot D_I
\end{equation}
where $\ii_I$ is the identity operator of $\hh_I$ and $\gamma^5$ the
graduation of the external geometry, is again a spectral triple. The
corresponding 1-forms are\cite{kt,vh}
$$
-i\gamma^\mu f_\mu^i \ot m_i + \gamma^5 h^j\ot n_j$$ where
$m_i\in\aa_I$, $h^j, f_\mu^i\in\cinf$, $n_j\in\Omega_I^1$.
Selfadjoint $1$-forms thus decompose in an  $\aa_I$-valued skew-adjoint
1-form field over $\mm$
\begin{equation}
\label{gaugg}
A_\mu\doteq f^i_\mu m_i,
\end{equation}
and an $\Omega^1_I$-valued selfadjoint scalar field 
\begin{equation}
\label{gaugh}
H\doteq h^j
n_j,
\end{equation}
 so that the part
of the covariant Dirac operator $D_A$ relevant for the distance formula is
\begin{equation}
\label{drelevant}
\mathcal{D} = D + \gamma^5
 H -i\gamma^\mu A_\mu  . 
\end{equation}

Assume now that $\aa_I$ is the algebra $\mn$ of $n$-square complex matrices. $\aa_E$ being nuclear the set of pure states of
\begin{equation}
\label{atot}
\aa =\cinf\ot \mn
= C^{\infty}(M,\mn)
\end{equation}
 is \cite{kadison}
 $$\pp(\aa) \simeq \pp(\aa_E) \times \pp(\aa_I)$$ where
 $\pp(\aa_I)\ni\omega_\xi$ is in $1$-to-$1$ correspondence with the projective
 space $\cc P^{n-1}\ni\xi$ via
\begin{equation}
 \label{evaluationun}
 \omega_\xi(m) = \scl{\xi}{m\xi} = \tr{s_\xi\, m}\quad \forall m\in\aa_I
 \end{equation}
with  $s_\xi\in\mn$ the support of $\omega_\xi$ (the density matrix in
physics) . The evaluation of
 \begin{equation}
\label{defxiox}
\xox \doteq
 (\omega_x\in\pp(\aa_E), \omega_\xi\in\pp(\aa_I))
\end{equation}
  on
 $a=f^i\ot m_i\in\aa$ reads
\begin{equation}
\label{evaluationdeux}
 \xi_x(a)= \tr{s_\xi\, a(x)}.
 \end{equation}
In other words $\pp(\aa)$ is a trivial fiber bundle 
$$P\overset{\pi}{\rightarrow}\nolinebreak M$$ 
with fiber $\cc P^{n-1}$ carrying a faithful action of $SU(n)$,
\begin{equation}
\label{freeaction}
\xox \mapsto {(U\xi)}_{x} \doteq \xox\circ\text{ad}(U^{-1})\quad\quad \forall U\in SU(n).
\end{equation} 

Back to (\ref{gaugg}),  $A_\mu$ has now value in the set $\frak{u}$
of skew-adjoint elements of of $\Pi(\aa_I)$.
Disregarding the distinction between the algebra and its
representation, 
\begin{equation}
\label{unim}
\frak{u} = \mathfrak{u}(n)
\end{equation}
 and $A_\mu$ can be
viewed as the local form of the
$1$-form field associated to some  Ehresmann connection $\Xi$ on the trivial
$U(n)$-principal bundle on $M$. Through a reduction to the trivial
$SU(n)$-principal bundle followed by a mapping\cite{kob} to the associated bundle $P$,
$\Xi$ yields a notion of horizontal curve on the bundle of pure
states, see (\ref{tangenthor}) below.
Note that in the rest of the
paper we focus on examples where $\Pi$ is the fundamental
representation of $\mn$ so that (\ref{unim}) is legitimate.

\subsection{Horizontal (or Carnot-Carath\'eodory) distance}

On the one side the gauge part $A_\mu$ of the covariant Dirac operator (\ref{drelevant})
equips the bundle $P$ of pure states with a horizontal
distribution. On the other side $A_\mu$ acquires a metric interpretation
 via the spectral distance $d$ (\ref{distancefluc}). A first
guess would be that the spectral distance coincides with
the natural distance associated to the horizontal distribution, namely
the Carnot-Carath\'eodory, or horizontal, distance $d_H$ defined
below. However as shown in [\citelow{cc}] and further studied in this
paper, the relation between $d$ and $d_H$
is more subtle. Before exploring this link in the following
subsection, we recall here some basic feature of the horizontal
distance. For interesting readers further deve\-lopments can be found in various texts on
subriemannian geometry, such as the excellent monography of Montgomery
[\citelow{montgomery}]. 

 Let $P\overset{\pi}{\rightarrow} M$ be a vector bundle over a
riemannian manifold with metric $g$. Assume $P$ comes equipped
with an Ehresmann connection $\Xi$, that is to say a way to split the tangent space into a
vertical subspace and a horizontal subspace,
\begin{equation}
\label{tangenthor}
T_pP = V_pP \oplus H_pP\quad\quad p\in P
\end{equation} where $HP$ is the kernel of the
connection $1$-form associated to $\Xi$. The horizontal (also called Carnot-Carath\'eodory) distance $d_H(p,q)$ is
defined as the length of shortest horizontal paths $t\in[0,1]\mapsto
c(t)\in P$ joining 
$p=c(0)$ to $q=c(1)$,
that is to say the length of the shortest path whose tangent vector is
everywhere horizontal with respect to the connection,
\begin{equation}
\label{dcc} d_H(p,q) \doteq \underset{\dot{c}(t)\in H_{\!c(t)}\!
P}{\text{Inf}}\; \int_0^1 \norm{\dot{c}(t)} dt \quad \forall
p,q\in P
\end{equation}
where the norm on $HP$ is the pull back of the metric
\begin{equation}
\label{normh}
\norm{\dot{c}} = \sqrt{g(\pi_*(\dot{c}),\pi_*(\dot{c}))}.
\end{equation}
Other choices of norm are possible, for instance one
could take the norm coming from
a metric defined on all $P$ (and not only on the base $M$). By the
choice of the pull-back (\ref{normh}) we actually define the length of
$c$ as the length of its projection $c_*$ on $M$. By convention when
$p,q$ cannot 
be reached by any horizontal path then
$d_H(p,q)$
is infinite.

\begin{figure}[hb]
\begin{center}
\hspace{-2truecm}
\mbox{\rotatebox{0}{\scalebox{.7
}{\includegraphics{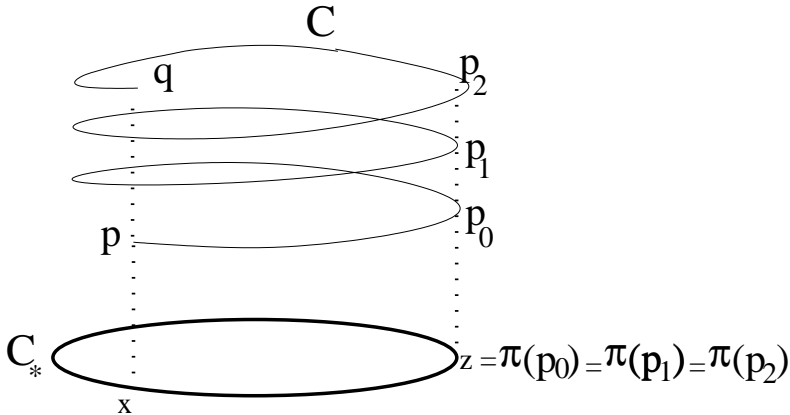}}}}
\end{center}
\caption{\label{helixfig} Horizontal lift of the circle.}
\end{figure}

As an illustration let us consider the trivial bundle ${\rr^2}^* \times \rr$ with 
local connection $1$-form
\begin{equation}
\label{exempleh}
A_\mu dx^\mu = (x^2 dx^1 - x^1dx^2) \ot \theta{\partial_3}
\end{equation}
 where $(x^1,x^2)$ are the coordinates of $\rr^2$, $\theta$ is a real
constant and $\partial_3$ is the base of the Lie algebra of the
additive group $\rr$. The corresponding connection $1$-form is
\begin{equation}
\label{exemplehbis}
\omega = A_\mu dx^\mu + \Theta 
\end{equation}
where $\Theta = dx^3 \ot\partial_3$ - $x^3$ the coordinate on $\rr$ -
denotes the Maurer-Cartan form. A curve $c(t) = (x^1(t),x^2(t),x^3(t))$ in the
bundle is horizontal if its tangent vector $\dot{c}$ is in
the kernel of $\omega$ that is, in polar coordinates $(r,\phi)$, if
\begin{equation}
 \label{contrainteh}
\dot{x_3} = \theta r^2\dot{\phi}.
\end{equation}
Therefore $p_0 = (z\in \rr^2, x^3)$ is linked to $p_1 =
 (z,x^3 + 2\abs{z}^2\theta\pi)$ by the horizontal lift $c$ of the
 circle $\mathcal{C}_z$ of radius $\abs{z}$ centered on the origin  (see figure \ref{helixfig}). By 
Lagrange multiplier one shows that under the constraint (\ref{contrainteh}) the length function
on the base
$$\int_{p_0}^{p_1} \sqrt{r^2\dot{\phi}^2 + \dot{r}^2}$$
has an extremum for $R = \text{constant}$. 
Assuming this extremum is a minimum (this point is discussed in the
next subsection), $c$ is the shortest horizontal
path between $p_0$ and $p_1$ and
$$d_H(p_0,p_1) = \text{perimeter of } \mathcal{C}_z = 2 \pi \abs{z}.$$

\subsection{Two distances from a single connection}

In the same way that the spectral distance on a manifold is
bounded by the geodesic distance, 
one can prove
[\citelow{cc}, prop III.1] that the gauge fluctuated spectral distance
(\ref{distancefluc}) in the product of a manifold by
$M_c(\cc)$ is bounded by the horizontal distance  
\begin{equation}
\label{dinfdh}
d(\xox, \yoz) \leq d_H(\xox, \yoz) \quad \forall \xox,\yoz\in P.
\end{equation}
In [\citelow{gravity}] was suggested that the two distances
were equal. However they are not because there do not
necessarily exists an equivalent of the function $L$ of
(\ref{functionl}) making $d_H$ the lowest upper bound for $d$ in
(\ref{dinfdh}). This has been shown in
great detail in [\citelow{cc}, prop. IV.4] and we briefly repeat here the
argument without proof. 

Identifying points to pure states,
eq.(\ref{functionl}) restricted to a minimal  geodesic between $x$ and
$z$ reads $$\omega_z(L) = d_\text{geo} (\omega_x,\omega_z).$$ In case of
a fiber bundle, one obtains instead that $d_H$ is the lowest upper bound for $d$
if and only if there exists some $a_0\in\aa$ such that
\begin{equation}
\label{holoun}
c_t(a_0) = d_H(\xox, c_t)
\end{equation}
for any $c_t\doteq c(t)$ on the minimal horizontal curve between $c(0)=\xox, c(1) =
\yoz$. Writing the explicit evaluation (\ref{evaluationdeux}) of $c_t$ on $a_0$
for each of the states $c_t = p_i$ that projects on the same point $z$ on the
base (fig. \ref{helixfig}), one obtains
\begin{equation}
\label{holonomycons}
p_i(a_0) = \tr s_{p_i} a_0(z) = d_H(\xox, p_i).
\end{equation}
Hence each point $p_i$ puts a condition on the same matrix $a_0(z)$ and
it is most likely that if there are too many such points, say more
than $n^2$, there is very little chance that an element $a_0$ satisfying
(\ref{holonomycons}) does exist.

\begin{figure}[hb]
\begin{center}
\hspace{-2truecm}
\mbox{\rotatebox{0}{\scalebox{1.3}{\includegraphics{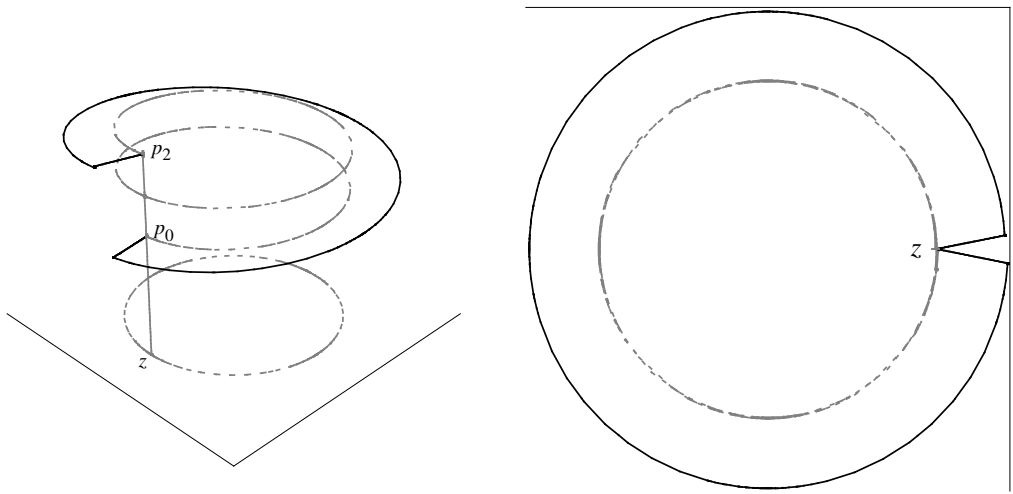}}}}
\end{center}
\caption{\label{doublevue} Two horizontal curves between $p_0$ and $p_2$
  and their projection on the base. Thin line is the double lift
  $c_1$ of the circle of radius $\abs{z}$. Thick line is $c_3$ and has
  only one pair of points that project onto the same point on the base.}
\end{figure}

 However one may argue that the number of
points $p_i$ could be reduced by choosing another minimal horizontal
curve with less points projecting onto the same given point of the base. For instance
getting back to the example $\rr^2\times \rr$ (\ref{exempleh}), there
exists at least two horizontal curve between $p_0$ and $p_2$. One, say $c_1$,
is the double lift of the circle $\mathcal{C}_z$ and has length
$$d_1 = 4\pi\abs{z}.$$ The other one is found by noting that
$x_3=\text{constant}, \phi=\text{constant}$
is a solution of the constraint (\ref{contrainteh}) so that the horizontal
lift (with initial condition
$x^3 =0$) of a radial segment $[z,\lambda z]$,
$\lambda\in\rr$, is the segment itself. Hence the horizontal lift  $c_2$  of
$$ [z,\sqrt{2} z]\circ {\cal C}_{\sqrt{2} z}\circ[\sqrt{2} z,z]$$
is also a horizontal path between $p_0$ and $p_2$ and has length
$$d_2 = 2\pi\sqrt{2}\abs{z} + 2\abs{z}(\sqrt{2}-1).$$ 
As $$d_1-d_2 = 2\abs{z}
(\sqrt{2}-1)(\sqrt{2}\pi -1)>0$$ 
$c_1$ is not the shortest horizontal path. Of course both $c_1$ and
$c_2$ have infinitely many pairs of points projecting onto the same
points on the base. But rotating a little bit the radial segment
around the $z$-vertical axis, it is not difficult to find a new
horizontal curve $c_3$ whose length is arbitrarily close to $d_2$ and
has no
such points (except $p_0$ and
$p_2$ themselves). Such a curve is plotted in figure \ref{doublevue}. 
  
Hence it is possible to deform a horizontal curve, keeping its
length almost fixed, in such a way to reduce the number of points
$p_i$ projecting onto the same $z$ in the base. However it is
difficult to extrapolate 
some general result when one has to take into account both the
connection and the metric on the base. Here we
found $c_2$  by moving from $z$ to
$\sqrt{2}z$, making a circle and come back to $z$. On a base with arbitrary
metric, it could be that one has interest to move far from $z$, make a
small circle there and come back. In fact even with the Euclidean metric a connection
$\frac{\omega}{r^2}$ instead of $\omega$ yields a different result: one has interest to go close to the origin, make a small
circle around it and come back. Moreover knowing there
exists a minimal horizontal curve whose projection does not self-intersects
is not sufficient to prove that the spectral distance coincide with
the horizontal distance. The holonomy obstruction (\ref{holonomycons})
is a necessary condition but certainly not a sufficient one to
guarantees that $d$ equal $d_H$. In the following we escape the
problem by considering a bundle on the circle,
$$M = S^1.$$
Thus, whatever the connection, there is at most one minimum
horizontal curve - a horizontal lift of $S^1$ with winding number $k$
- between any two
points $\xox, \xoz$ on the same fiber. Then 
\begin{equation}
\label{dhn2}
d_H(\xox,\xoz) = 2k\pi
\end{equation}
and computing explicitly $d(\xox,\xoz)$ one finds that both distances
do not coincide. This has been shown in [\citelow{cc}] for low dimension 
($\aa_I = M_2(\cc)$) with constant connection
($A_\mu = diag(0,-i\theta\in i\rr)$). In the rest of the paper
we generalize this result to non constant connections of dimension $n\geq 2$.

 \section{Connected components\label{composanteconnes}}   

In this section we work out the topological aspect of both the spectral and
the horizontal distances induced by the fluctuation of the metric in a
spectral triple built on $\aa = C^\infty(S^1,\mn)$. Main results are
propositions \ref{proptauxi} and \ref{propconneccomp} that are further
discussed in section \ref{twosection}.

\subsection{A noncommutative geometry for the circle}

Consider the product (\ref{tripletprod}) of the natural spectral
triple (\ref{continu}) on $M = S^1$ by the finite dimensional spectral triple
$$\aa_I = M_n(\cc),\quad \hh_I = M_n(\cc),\quad D_I = 0.$$
Note that $\aa_I$  is represented on itself - by simple matrix
multiplication - rather than on $\cc^n$ in order to fulfill the reality axiom of spectral
triple\cite{gravity}. However the later does not enter the
distance formula and the operator norm of
an $n$-square 
matrix acting on $\mn$ coincides with its operator norm on $\cc^n$. Consequently
regarding distance computation it is equivalent to consider that
$\aa_I$ acts on $\cc^n$.  The vanishing of $D_I$
guarantees that the scalar part $H$ of the fluctuation is zero, so
that the part of the fluctuated Dirac operator (\ref{drelevant})
relevant for the spectral distance is the usual covariant operator
\begin{equation}
\label{tildedeux}
 {\cal{D}} = -i\gamma^\mu (\partial_\mu \ot \ii_I  + \ii_E\ot A_\mu)
\end{equation}
associated to the connection $A_\mu$ on the product bundle
$L_2(M,S)\ot(M\times\mn) = L_2(M,S\ot \mn)$. Strictly speaking one cannot talk
 of a covariant operator on the associated $SU(n)$-bundle
$$\pp(\aa) = P\overset{\pi}{
\rightarrow} S^1$$ of pure states of $\aa$
 since the fiber $\cc P ^{n-1}$ is not a vector
space.
However as explained below eq.(\ref{unim})  $A_\mu$ does define a
notion of horizontal curve in $P$. Technical details on that matter
are given in the next subsection. 

Two states in $P$ are at infinite horizontal distance if and only if
there is no horizontal path between them. All pure states at finite
horizontal distance from $\xox$ are said {\it accessible},
\begin{equation}
\label{acc} \text{Acc}(\xox)\doteq \{q\in P;\; d_H(\xox,q) <
+\infty\}.
\end{equation}
Two pure states $\xox$, $\yoz$ are said connected for $d$ if and
only if $d(\xox, \yoz)$ is finite and we write
\begin{equation}
\label{con} \text{Con}(\xox)\doteq \{q\in P;\; d(\xox,q) <
+\infty\}
\end{equation}
the connected component of $\xox$ for $d$. From (\ref{dinfdh})
$\text{Acc}(\xox)$ is connected for $d$,
\begin{equation}
\label{inclusion} \text{Acc}(\xox)\subset\text{Con}(\xox).
\end{equation}
 However there is no reason for $\text{Con}(\xox)$ to equal
$\text{Acc}(\xox)$. $d$ can remain finite although $d_H$ is
infinite. As discussed in the precedent section  what matters is the holonomy of the connection.
In low dimension ($n = 2$) with constant connection
 we showed in [\citelow{cc}]
that $\text{Con}(\xox)$ is a $2$-torus and $\text{Acc}(\xox)$ a
dense subset if $\theta$ is irrational. In the rest of this section
we show that for non-constant connection in arbitrary $n$ dimension 
$\text{Con}(\xox)$ is an $n_c$ torus $\mathbb{U}_\xi$ with $n_c\leq n$ given by the
holonomy. $\text{Acc}(\xox)$ is a dense, discrete or a discrete union of dense
subsets of $\mathbb{U}_\xi$.

\subsection{Accessible points}
Let us begin by $\text{Acc}(\xox)$, which is nothing but the horizontal lift of the circle with initial condition given by $\xi$.
Here ``horizontal lifts'' take place in the associated bundle $P$ of pure
states. So, as a starting point, it may be useful to recall how the horizontal structure defined on the trivial $SU(n)$-bundle by the
gauge fluctuation $A_\mu$ is exported to $P$. Explicitly a
curve in $P$ is horizontal if\cite{kob} it writes as the product $v U$ of

- a curve $v$ in $P$  whose projection on the base
coincides with the projection ${c_*}$ of $c$. 

- a curve $U$ in $SU(n)$ such that $U^{-1}\dot{U} =
-\omega(\dot{v})$ where $\omega$ is the connection $1$-form defined
by $A_\mu$ via (\ref{exemplehbis}). Locally these write
\begin{equation}
\label{horizontalrig}
A_\mu \dot{c}_*^\mu = - U^{-1}\dot{U}
\end{equation}
and the horizontal lift $c\subset P$ of a given $c_*\in M$ with initial
conditions $c_*(0)=x,\, c(0) = \xi_x$ writes 
\begin{equation}
\label{releverig}
c(t) = (U_t\xi)_y
\end{equation}
where $y=c_*(t)$, $U_t$ is obtained by integrating
(\ref{horizontalrig}) and the action of $SU(n)$ on $\cc P^{n-1}$ is given in
(\ref{freeaction}). In our specific case of interest $M=S^1$ the connection $1$-form has only one component
\begin{equation}
\label{amunique}
A_\mu = -A_\mu^* = A. 
\end{equation}
Once for all we fix on $\hh_I$
a basis of real eigenvectors of $iA$ such that
\begin{equation}
\label{diagA}
A = i\dm{ccc} \theta_{1} &\ldots & 0\\  \vdots & \ddots&\vdots\\
0&\ldots & \theta_n\fm \end{equation} where the $\theta_j$'s are real
functions on $S^1$. Let $\rr \text{ mod } 2\pi$ parametrize the circle and
$x$ be the point with coordinate $0$. Any curve $c_*$ in $S^1$ starting at
$x$ is 
$$c_*(t) = t \quad\forall t\in [0,\tau],\quad \tau\in\rr$$ where
we do not restrict to $\tau\leq 2\pi$ since we want to include curves
with winding numbers greater than $1$.
Then
\begin{equation}
\label{orderedpath}
U_\tau = \text{exp}\left(-\int_0^\tau A(t) dt \right)
\end{equation} 
acts on $\xi\in\cc P^{n-1}$ - viewed as an $n$-complex normalized
vector - by usual matrix multiplication.
 Within a trivialization $(\pi,
V)$ the initial condition writes
 $$V(c(0))=\xi = \dm{c} V_1\\\vdots \\ V_n\fm \in \cc P^{n-1}$$
 and the horizontal lift $c(\tau)=(c_*(\tau), V(\tau))$ has components
\begin{equation}
\label{vi} \displaystyle V_j(\tau) = V_j e^{-i\Theta_j(\tau)}
\end{equation} with
\begin{equation} \label{gdtheta} \Theta_j(\tau) \doteq \int_0^ {\tau}
\theta_j(t)dt.
\end{equation}

\begin{figure}[bht]
\begin{center}
\mbox{\rotatebox{0}{\scalebox{.6}{\includegraphics{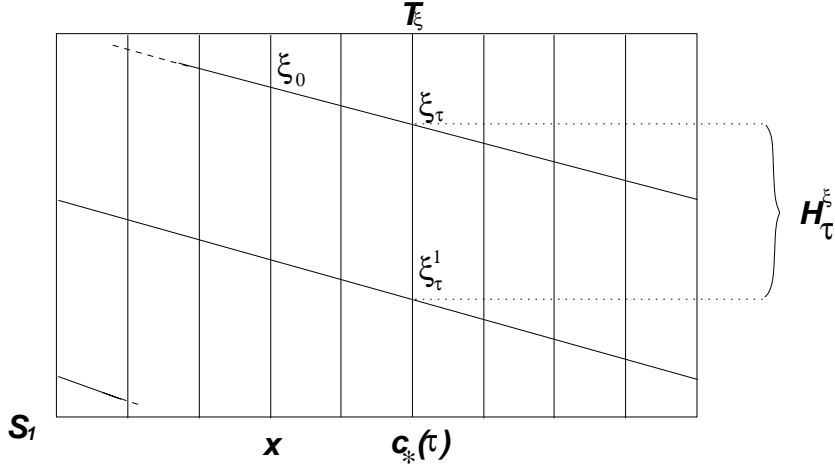}}}}
\caption{{The $n$-torus $\mathbb{T}_\xi$ with base $S^1$ and fiber the
    $(n-1)$-torus $T_\xi$. The diagonal 
line is $\text{Acc}(\xox =\xi_0)=\left\{\xi_\tau,
  \tau\in\rr\right\}$.  $H_\tau^\xi$ is the orbit of $\xi_\tau$
    under the holonomy group at $c_*(\tau)$.}}
\label{figtxi}
\end{center}
\end{figure}

In other terms, with notation (\ref{defxiox}), the elements of $P$ accessible from 
$\xi_x = \xi_0$
 are the pure states
\begin{equation}
\label{xitau} \xi_\tau\doteq (\omega_{c_*(\tau)},
\omega_{V(\tau)}),\quad \tau\in\rr.
\end{equation}
On a given fiber $\pi^{-1}(c_*(\tau))$, $\text{Acc}(\xox)$ reduces
to 
\begin{equation}
\label{xitauk}
\xi_\tau^k\doteq \xi_{\tau + 2k\pi},\quad k\in\zz,
\end{equation}
with components
\begin{equation}
\label{xitaukun}
\displaystyle
V_j(\tau + 2k\pi) = V_j(\tau) e^{-ik\Theta_j(2\pi)}.
\end{equation} Dividing each
$\xi_\tau^k$ by the irrelevant phase $e^{-ik\Theta_1(2\pi)}$ and
writing
\begin{equation}
\label{thetaij} \Theta_{ij} \doteq \Theta_i -\Theta_j,
\end{equation}
 one obtains the orbit $H_\tau^\xi$ of $\xi_\tau$ under the action of the
holonomy group at $c_*(\tau)$, namely
$$H^\xi_\tau \doteq \text{Acc}(\xox)\cap \pi^{-1}(c_*(\tau)) =
\{ \dm{r} V_1(\tau)\\
e^{ik\Theta_{1j}(2\pi)}V_j(\tau) \fm,\; k\in\zz; j=2,..., n\}.$$
Fiber-wise the set of accessible points is thus a subset of the
$(n\!-\!1)$-torus of $\cc P^{n-1}$,
\begin{equation} \label{txi}
 T_\xi \doteq  \{ \dm{r} V_1\\
e^{i\varphi_j}V_j \fm,\; \varphi_j\in \rr,\; j=2,..., n\},
\end{equation}
where we noted that $V_i(\tau)$ differs from $V_i$ only by a phase
factor. The union on all fibers yields:
\begin{proposition}
\label{proptauxi}
$\text{Acc}(\xox) = \underset{\tau\in [0,2\pi[}{\bigcup} H_\tau^\xi$ is a
subset of the $n$-torus (see figure \ref{figtxi})
\begin{equation}
\label{tz0}\mathbb{T}_\xi\doteq S^1 \times T_{\xi}.
\end{equation}
\end{proposition}

\begin{figure}[hb]
\hspace{-.9truecm}
\mbox{\rotatebox{0}{\scalebox{.8}{\includegraphics{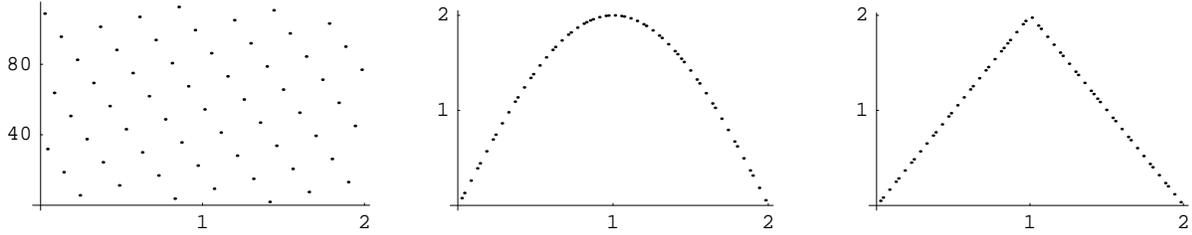}}}}
\caption{ \label{dss} from left to right, the horizontal, spectral and
  euclidean distances on the fiber in the $n=2$ case.
The fiber $T_\xi = S^1$ is
parametrized by $\varphi\in [0,2\pi]$ and the plots is the distance
between $\varphi=0$ and $varphi = k\Theta_{12}(2\pi$). 
}
\end{figure}

Let us comment on this inclusion 
$$\text{Acc}(\xox)\subset \mathbb{T}_\xi\subset P$$
which fiber-wise reads
$$H_\tau^\xi\subset T_\xi\subset \cc P^{n-1}.$$ 
 Two elements of $\cc P^{n-1}$ with components $V_i$'s, $U_i$'s belong to different tori $T_\xi$ if
and only if $\abs{V_i}\neq \abs{U_i}$ for at least one $i\in [2,n]$.
In this sense the $T_\xi$'s realize a slicing of $\cc P^{n-1}$, and  the
$\mathbb{T}_\xi$'s a slicing of $P$, indexed
by the modules of the components of $\xi$. The latter does
not correspond to the horizontal foliation of $P$
for
\begin{equation}
\label{acctxidiff} \text{Acc}(\xox)\varsubsetneq \mathbb{T}_\xi.
\end{equation} 
 Indeed at best, when all the $\Theta_{1j}(2\pi)$'s are distinct and
irrational, $\mathbb{T}_{\xi}$ is the completion of
${\text{Acc}(\xi_x)}$ with respect to the euclidean norm on $T_\xi$.
But as soon as one of the $\Theta_{1j}(2\pi)$'s is rational,
$\text{Acc}(\xox)$ is no longer dense in $\mathbb{T}_\xi$. For
instance when all the $\Theta_{1j}(2\pi)$'s are rational,
$H^\xi_\tau$ is a discrete subset of $T_\xi$. In any case $\text{Acc}(\xox)$ cannot be the whole
torus $\mathbb{T}_{\xi}$. Later on the inclusion (\ref{acctxidiff})
will be refined by the computation of the connected component for the
spectral distance. For the moment let us simply underline that the
horizontal distance "forgets" about the fiber bundle structure of
the set of pure states; and it forgets it twice: 
\begin{enumerate}
\item  two pure states on the same fiber belonging
to distinct tori $\mathbb{T}_\xi$, $\mathbb{T}_\eta$ are infinitely Carnot-Carath\'eodory far
from each other although they might be close in a suitable topology
of $\cc P^{n-1}$.\newpage 
\item  on a given component $\mathbb{T}_\xi$, 
$d_H(\xi_\tau, \xi_\tau^k)  = 2k\pi$.
So in case the $\Theta_{1j}(2\pi)$'s are
irrational, by (\ref{xitauk}) and (\ref{xitaukun}) one can find
close to $\xi_\tau$ (in the euclidean topology of $\mathbb{T}_\xi$)
some $\xi_\tau^k$ arbitrarily Carnot-Carath\'eodory far from
$\xi_\tau$ (see fig. \ref{dss}).
\end{enumerate}
On the contrary the spectral distance keeps better in
mind the fibern structure of $P$. For instance in the case $n=2$
(i.e. $T_\xi = S^1$) with a constant non-trivial connection (i.e.
$\theta_1\neq\theta_2\in\cc$) we found in [\citelow{cc}] that
$\text{Con}(\xox)$ is the whole $2$-torus $\mathbb{T}_\xi$. Two
states on different tori $\mathbb{T}_\xi$, $\mathbb{T}_\zeta$ are
still infinitely far from one another, but on a given $\mathbb{T}_\xi$ all
states are at finite spectral distance from one another.
Specifically on each fiber the spectral distance $d$
appears as a "smoothing"\cite{smoother} of the euclidean distance $d_E$ on the
circle (figure \ref{dss}).
\vspace{.5cm}

\subsection{Connected points} 

Let us investigate now the connected component $\text{Con}(\xox)$ of the
spectral distance in an $\mn$ bundle over the circle with arbitrary
connection. As
shown in proposition \ref{propconneccomp} below, $\text{Con}(\xox)$
is a subset of the torus $\mathbb{T}_\xi$ defined in (\ref{tz0}). But
contrary to the $n=2$ case\cite{cc} $\text{Con}(\xox)$ is not necessarily equal to
$\mathbb{T}_\xi$. Viewing the torus $\mathbb{T}_\xi$ as the subset
of $\rr^n$, 
\begin{equation}
\mathbb{T}_\xi = \left\{ \tau\in [0,2\pi[,\varphi_i \in [0,2\pi[, i=2,...n\right\}
\end{equation} 
we show below that  $\text{Con}(\xox)$ is a sub-torus
$\mathbb{U}_\xi$ of $\mathbb{T}_\xi$,
\begin{equation}
\mathbb{U}_\xi = \left\{ \tau\in [0,2\pi[,\varphi_i \in [0,2\pi[,i=2,...n_c\right\}
\end{equation}
with dimension $n_c\leq n$ given by the number of equivalence classes of the following
relation. 
 
\begin{definition}
\label{connectdirect} We say that two directions $i,j$ of
$\mathbb{T}_\xi$ are {\it far} from each other if the components $i$
and $j$ of the holonomy at $x$ are equal, and we write
$\text{Far}(.)$ the equivalence classes,
\begin{equation}
\text{Far}(i)\doteq \{ j\in [1,n] \text{ such that } \Theta_j(2\pi)
= \Theta_i(2\pi)\;\text{mod}[2\pi]\}.
\end{equation}
Two directions belonging to distinct equivalence classes are said
close to each other. We denote $n_c$ the numbers of such equivalence
classes and we label them as
$$\text{Far}_1 = \text{Far}(1),\,
\text{Far}_p = \text{Far}(j_p)\quad p=2,... ,n_c$$ where $j_p\neq0$
is the smallest integer that does not belong to
$\underset{q=1}{\overset{p-1}{\bigcup}}\text{Far}_q$.
\end{definition}

\noindent Our terminology will become clearer after proposition
\ref{propconneccomp}. The word ``far'' simply reflects that some directions (one may prefer say
``dimensions'') of $\mathbb{T}_\xi$ will be found to be infinitely
\emph{far} from one another
with respect to the spectral distance. Identifying them by the equivalence relation
III.2, we will see that each equivalence class contribute to
$\text{Con}(\xox)$ by furnishing one of the dimension of
$\mathbb{U}_\xi$. This is explained in section \ref{twosection} below and
illustrated in figure \ref{troistore}.

 Before establishing this result, let us recall two
simple lemmas: first [\citelow{finite}, lemma 1] the supremum in the
distance formula can be searched on selfadjoint elements of $\aa$.
Second [\citelow{cc}, lemma III.2]
\begin{lemma}
 \label{infini}
 $d(\xox, \yoz)$ is infinite if and only if
there is a sequence $a_n\in\aa$ such that
\begin{equation}
\label{condinf} \underset{n\rightarrow +\infty}{\text{lim}}
\norm{[\dd,a_n]} = 0,\quad \underset{n\rightarrow
+\infty}{\text{lim}} \abs{\xox(a_n)-\yoz(a_n)} = +\infty.
\end{equation}
\end{lemma}

\noindent This allows to prove the main result of this section:

\begin{proposition}\label{propconneccomp}
$\text{Con}(\xox)$ is the $n_c$ torus
\begin{equation}
\label{defu} \mathbb{U}_\xi\doteq \underset{\tau\in [0,2\pi[}{\bigcup}
U^\xi_\tau \end{equation} where $U^\xi_\tau\subset T_\xi$ is the
$(n_c\!-\!1)$ torus defined by ($V_{i}(\tau)$ is given in
(\ref{vi}))
\begin{equation}
\label{txif} U^\xi_\tau \doteq  \{ \left( \begin{array}{rl}
 V_{i}(\tau) & \forall
i\in\text{Far}_1\\
e^{i\varphi_2}V_{i}(\tau) & \forall i\in\text{Far}_2\\
\ldots & \\
e^{i\varphi_{n_{\!c}}}V_{i}(\tau) & \forall i\in\text{Far}_{n_c} \fm
, \varphi_j\in\rr, j\in [2,n_c]\}.
\end{equation}
\end{proposition}

\noindent\noindent{\it Proof.} The idea is to determine by
lemma above all states that are at infinite spectral distance from
$\xox$. To do so we first need to identify the $a_0\in\aa$ that
commute with $\mathcal{D}$. Let $a_{ij}\in \aa_E$ be the
components of $a=a^*\in\aa$, identified as $2\pi$-periodic complex
functions on $\rr$,
\begin{equation}
\label{aij}
 a_{ij}(\tau)\doteq a_{ij}(c_*(\tau)) = a_{ij}(\tau +
2k\pi)\quad k\in\zz \end{equation} with initial condition
$a_{ij}(0)= a_{ij}(x)$. Let dot denote the derivative. The
Clifford action reduces to the multiplication by $1$ ($\gamma^\mu
= \gamma^1 = 1$) so that $[D_E, a_{ij}] = -i\dot{a_{ij}}$ and
\begin{equation}
\label{diraac} i[\dd,a] =  \dm{cccc}
\dot{a_{11}}&\dot{a_{12}} +ia_{12}\theta_{12}&\ldots& \dot{a_{1n}} + ia_{1n}\theta_{1n} \\
\dot{a_{21}}  - ia_{21}\theta_{12}& \dot{a_{22}}& &\vdots\\
\vdots & & \ddots& \\\dot{a_{n1}} - ia_{n1}\theta_{1n} &\ldots
&\ldots &\dot{a_{nn}} \fm
\end{equation}
where $\theta_{ij}\doteq \theta_i -\theta_j$. The
commutator (\ref{diraac}) is zero if and only if
\begin{equation}
\label{comun} a_{ii} = C_i = \text{constant, }\;
a_{ij}(\tau) = C_{ij}e^{-i\Theta_{ij}(\tau)}\end{equation} 
where $$\Theta_{ij}(\tau) = \int_0^\tau \theta_{ij}(t)\,dt$$ 
has already been
defined in (\ref{gdtheta},\ref{thetaij}). 
The $a_{ij}$'s being $2\pi$-periodic functions (\ref{comun})
has non zero solutions $C_{ij}$ only if
\begin{equation}
\Theta_{ij}(\tau + 2k\pi) = \Theta_{ij}(\tau) \,\text{mod}[2\pi]
\end{equation}
that is to say, since by periodicity of $\theta_{ij}$
\begin{equation}
\label{lineairetheta} \Theta_{ij}(\tau + 2k\pi) =  \Theta_{ij}(\tau)
+ k\Theta_{ij}(2\pi), \end{equation} only if the directions $i$ and
$j$ are far from each other in the sense of definition
\ref{connectdirect}. In other terms $a_0$ commutes with $\dd$ if and
only if its components $a_{ij}$ vanishes for any couples
of close directions $i,j$ and are given by (\ref{comun}) otherwise.

For such $a_0$ and for any $\xi\in\cc P^{n-1}$ with components $V^i$,
(\ref{evaluationdeux}) and (\ref{evaluationun}) 
with $y=c_*(\tau)$ yield
\begin{equation} \label{explicite} \yox(a_0) =
\overset{n}{ \underset{i=1}{\sum}}C_i\abs{V_i}^2 +
\underset{i<j}{\sum} 2\re \lp \bar{V_i}V_j
C_{ij}e^{-i\Theta_{ij}(\tau)}\rp.
\end{equation}
Consider now $\zeta\in\cc P^{n-1}$ with components $W_i$ such that,
for at least one value $i_0\in[1,n]$,
\begin{equation} \label{casun}
\abs{V_{i_0}}\neq\abs{W_{i_0}}. \end{equation} One easily finds some
$a_0$ such that $\xox(a_0) \neq \yoz(a_0)$, for
instance
$$C_{ij}=C_i = 0  \text{ for all } i,j \text{ except } C_{i_0}\neq
0.$$ Thus $d(\xox, \yoz)$ is infinite by lemma \ref{infini}
(consider $a_n = na_0$). This means that
$\text{Con}(\xox)\cap\pi^{-1}(c_*(\tau))$ is included within the set
of pure states $\yoz$ that do not satisfy (\ref{casun}), that is to
say those
with components
\begin{equation}
\label{casdeux}
 W_i = e^{i\varphi_i}V_{i},\quad \varphi_i\in\rr,\, i=1, ...,n.
\end{equation}
Dividing by an irrelevant phase $e^{i\varphi_1}$, one is
back to the torus $T_\xi$ defined in (\ref{txi}). This indicates that
$\text{Con}(\xox) \subset \mathbb{T}_\xi$ as expected. 

However there is no
reason for $\text{Con}(\xox)$ to be the whole of
$\mathbb{T}_\xi$. Intuitively this is because the sum in
(\ref{explicite})
runs only on couples of directions that are far from each other. The
vanishing of $a_{ij}$ for close directions indicates that the latter cannot
contribute to some $a_0$ making $d(\xox, \yoz)$ infinite, i.e. such that $\xox(a_0) \neq  \yoz(a_0).$ To make
this idea more precise let us first observe  that $V_i$
differs from $V_i(\tau)$ only by a phase factor. Therefore
redefining
\begin{equation}
\label{phiredef}
\varphi_i \rightarrow \varphi_i - {i\Theta_i(\tau)}
\end{equation}
 condition (\ref{casdeux})  equivalently writes
\begin{equation}
\label{defw}
 W_{i} = e^{i\varphi_i}V_{i}(\tau).
\end{equation}
Assuming $\zeta$ has such components (\ref{explicite}) yields
\begin{equation}
\label{diffun} \yoz(a_0) -\xox(a_0) = \underset{i<j, j\in\text{Far}(i)}{\sum} 2\re
\lp\bar{V_i}V_j C_{ij}( e^{-i(\varphi_i - \varphi_j)}-1)\rp
\end{equation}
where $e^{-i\Theta_{ij}(\tau)}$ has been absorbed in the redefinition
(\ref{phiredef}) and  $1 = e^{-i\Theta_{ij}(0)}$.
  By lemma \ref{infini} $d(\xox,\yoz) = +\infty$ if and only if there
  exists at least one element $a_0$ (or a sequence $a_n$ such that
  $\underset{n\rightarrow +\infty}{\text{lim}} {{a_{n}}_{ij}}$ has
  the same properties of ${a_0}_{ij}$) such that (the limit of)
  (\ref{diffun}) 
does not vanishes. This happens if and only if 
$$\varphi_i \neq \varphi_j  \,\text{mod}[2\pi]$$ 
for at least one $j$ in $\text{Far}(i).$ Thus $\text{Con}(\xox)\cap\pi^{-1}(c_*(\tau))$
is the set of states that do not satisfy condition above, i.e. those with components
$$
W_{i} = e^{i\varphi_j}V_{i}(\tau) \quad\, \forall
i\in\text{Far}(j).$$ 
Up to an irrelevant phase
$e^{i\varphi_1}$, this is nothing but $U^\xi_\tau.$\hfill $\blacksquare$
\newline

\subsection{\label{twosection}Two topologies from a single connection}

\begin{figure}[hbt]
\begin{center}
\mbox{\rotatebox{0}{\scalebox{.6}{\includegraphics{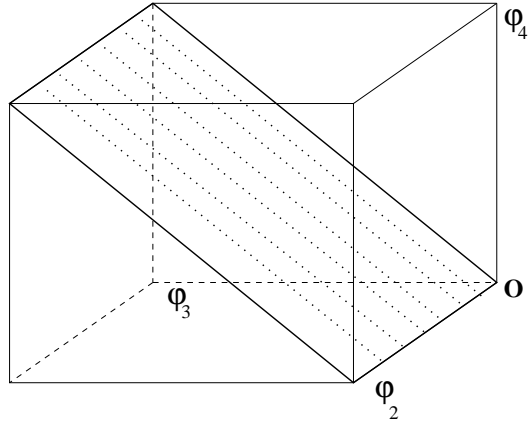}}}}
\end{center}
\caption{\label{troistore} The connected components for the spectral
  and the horizontal distances on the fiber over $x$, in case $n=4$ with directions
  $3,4$ far from each other.  We chose a rational $\Theta_{12}(2\pi)$
 and irrational $\Theta_{13}(2\pi) = \Theta_{14}(2\pi)$.
 $\abs{V_i}$ determine a $3$-torus
  $T_\xi=\left\{\varphi_i\in [0,2\pi[ , i=2,3,4 \right\} \subset\cc
  P^{3}$. $\text{Arg}\, V_i$ 
determine a point $O$ in $T_\xi$.
\newline
-$\text{Con}(\xox)\cap\pi^{-1}(x)$ is the $2$-torus $U_\xi = \left\{\varphi_2\in
    [0,2\pi [, \varphi_3=\varphi_4\in [0,2\pi [\right\}$ containing $O$.\newline
  -$\text{Acc}(\xox)\cap\pi^{-1}(x)$ 
is a subset of $U_\xi$, discrete in $\varphi_2$ and
dense in $\varphi_3\sim \varphi_4$.}
\end{figure}

The spectral and the horizontal distance yield two distinct topologies
$\text{Con}$ and $\text{Acc}$ on the bundle of pure states $P$. Obviously
$e^{i\Theta_{1j}(2k\pi)} = e^{i\Theta_{1i}(2k\pi)}\, \forall j\in\text{Far}(i)$, hence $H_\tau^\xi\subset U_\tau^\xi$ fiber-wise and
$\text{Acc} (\xox) \subset \mathbb{U}(\xox)$ globally, as
expected from (\ref{inclusion}). Also obvious is the inclusion of
$\mathbb{U}_\xi$ within $\mathbb{T}_\xi$. To summarize the various
connected components organize as follows,
\begin{equation}
\label{composantesconnexes} \text{Acc}(\xox)\subset
\text{Con}(\xox) = {\mathbb{U}}_\xi\subset {\mathbb{T}}_\xi\subset P
\end{equation}
or fiber-wise
\begin{equation}
\label{composantesconnexes} H_\tau^\xi \subset U_\tau^\xi\subset
T_\xi\subset \cc P^{n-1}.
\end{equation}
The difference between $\text{Acc}(\xox)$ and ${\mathbb{U}}_\xi$ is
governed by the irrationality of the connection - as explained below
eq. (\ref{acctxidiff}) - whereas the difference between
 ${\mathbb{U}}_\xi$ and $\mathbb{T}_\xi$ is governed by the number of
close directions. 
 More specifically
\begin{equation}
\label{diverstxi} \mathbb{T}_\xi= \underset{\zeta\in
T_\xi}{\bigcup}\text{Acc}(\zeta_x) \end{equation}
 is the union of all states with equal components up to phase
 factors. Meanwhile
\begin{equation}
\label{diversuxi} \mathbb{U}_\xi= \underset{\zeta\in
U_\xi}{\bigcup}\text{Acc}(\zeta_x), \end{equation}
 with
$U_\xi = U_{\tau = 0}^\xi$, is the union of all states with equal components up
to phase factors, with the extra-condition that phase factors
corresponding to directions far from each other must be equal. This is
illustrated fiber-wise in low dimension in figure \ref{troistore}.
 When 
all the directions are close to each other
(e.g. when the functions $\theta_i$'s are constant and distinct from
one another) then $U_\xi=T_\xi$, $n_c =n$ and
$$\mathbb{U}_\xi= \mathbb{T}_\xi.$$
On the contrary when all the directions are far from each other,
that is to say when the holonomy is trivial, $n_c = 1$ and, in
agreement with [\citelow{cc}, Prop. IV.1],
$${\mathbb{U}}_\xi =  S^1 = \text{Acc}(\xox).$$
In short $\text{Con}(\xox)$ varies from
$\text{Acc}(\xox)$ to $\mathbb{T}_\xi$ while 
$\text{Acc}(\xox)$ varies from a discrete to a dense subset of  $\text{Con}(\xox)$.

Note that none of the distance is able to "see" between different
tori $\mathbb{T}_\xi$, $\mathbb{T}_\eta$. However within a given
$\mathbb{U}_\xi$ the spectral
distance sees between the horizontal components.
In this sense the spectral distance keeps "better in mind"
the bundle structure of the set of pure states $P$.  This suggests
that the spectral 
distance could be relevant to study some transverse metric structure in a more general
framework of foliation.

\section{Spectral distance on the circle with a $\cc P^1$ fiber}

In the rest of the paper we investigate the spectral distance $d$ on
$\text{Con}(\xox)$ for arbitrary connection and any choice of
$\xox$. In the present section we focus on the case $n=2$ and explicitly compute $d$ on all
$\text{Con}(\xox)$, thus generalizing the result of [\citelow{cc}]
to non-constant connection and non-equatorial $\xi$ (see definition
below, eq.(\ref{equatorial})). In section 5 we compute $d$ on the
fiber for dimension  $n\geq 2$.

\subsection{Parametrization of $\mathbb{T}_\xi$}

The set of pure states of $\aa_I=M_2(\cc)$ identifies to the sphere
$S^2$ via the traditional correspondence
\begin{equation}
\label{defz} \xi=\lp\begin{array}{c} V_1\\
V_2\end{array}\rp \in\cc P^{1} \longleftrightarrow
\left(\begin{array}{l} x_\xi \doteq 2\Re (V_1\overline{V}_2)\\
y_\xi \doteq 2\Im(V_1\overline{V}_2)\\ z_\xi \doteq \abs{V_1}^2 -
\abs{V_2}^2\end{array}\right).
\end{equation}
The torus $T_\xi$ defined in (\ref{txi}) is now a circle $S_R$ of
radius $R\doteq (1-z_\xi^2)^{\frac{1}{2}}$ located inside $S^2$ at
the "altitude" $z_\xi$.\footnote{$T_\xi$ is characterized by a
  single parameter because the value of $z_\xi$ is enough to determine
$\abs{V_1}$ and $\abs{V_2}$.} The bundle of pure states of
$C^\infty(S^1,\m2)$ is sliced by $2$-tori
$$\mathbb{T}_\xi=S^1 \times S_R$$
each of them containing the connected component $\text{Con}(\xox) = \mathbb{U}_\xi$
 associated to the pure states
 of altitude $z_\xi$. 
This connected component is 
particularly simple: either the connection $A$ is a multiple of the
 identity and $\mathbb{U}_\xi = S^1$, or  $A$ is non-trivial so that the two directions of $\mathbb{T}_\xi$ are close
to each other in the sense of definition \ref{connectdirect} and
$\mathbb{U}_\xi = \mathbb{T}_\xi$ by proposition \ref{propconneccomp}. In [\citelow{cc}]
we computed $d(\xox,\yoz)$  for any
$\yoz\in\mathbb{T}_\xi$, assuming the connection was constant and $\xi$
 was an equatorial state, namely
 \begin{equation}
\label{equatorial}
 z_\xi =0.
 \end{equation}
In the following we investigate the spectral distance between any two
states in $\mathbb{T}_\xi$, equatorial or not, with arbitrary connection. 

A preliminary step is to suitably parametrize
$\mathbb{T}_\xi$. Recall that the base $S^1$ is parametrized by
$[0,2\pi[$ with $x$ the point of coordinate $0$. Any pure state accessible from $\xox$
 is $\xi_\tau$, $\tau\in\rr$, with
components $V_i(\tau)$ given in (\ref{vi}). In the fiber containing $\xi_\tau$, the $1$-torus
$T_\xi$ is obtained multiplying $V_2(\tau)$  by an
arbitrary phase $e^{i\varphi}$. Remembering (\ref{xitaukun}), one
obtains the same element of $T_\xi$  multiplying the components $V_1$,
$V_2$ at $\tau+2k\pi$ by 
 $e^{i k\Theta_1(2\pi)}$, $e^{i(\varphi +  k\Theta_2(2\pi))}$ respectively.
Equivalently one can let $V_1$ unchanged and
multiply $V_2(\tau + 2k\pi)$ by $e^{i(\varphi +
  k\Theta_{21}(2\pi)}$. Hence the following parametrization of $\mathbb{T}_\xi$.
{\definition{Given $\xox$ in $P$, any pure state $\yoz$ in the
    $2$-torus $\mathbb{T}_\xi$ is in
one-to-one correspondence with an equivalence class
\begin{equation}
\label{parameters} (k\in\nn,\; 0\leq \tau_0\leq 2\pi,\;0\leq\varphi
\leq 2\pi)\sim (k + \zz,\; \tau_0,\;\varphi- 2\zz\omega\pi)
\end{equation} with
\begin{equation}
\label{tauk} \tau = 2k\pi + \tau_0,\quad \omega\doteq
\frac{\Theta_1(2\pi) - \Theta_2(2\pi)}{2\pi}
\end{equation}
such that
\begin{equation}
\label{zetay}
\yoz =
\dm{r} V_1(\tau)\\
e^{i\varphi}V_2(\tau)\fm.
\end{equation}}}

\begin{center}
\mbox{\rotatebox{0}{\scalebox{.5}{\includegraphics{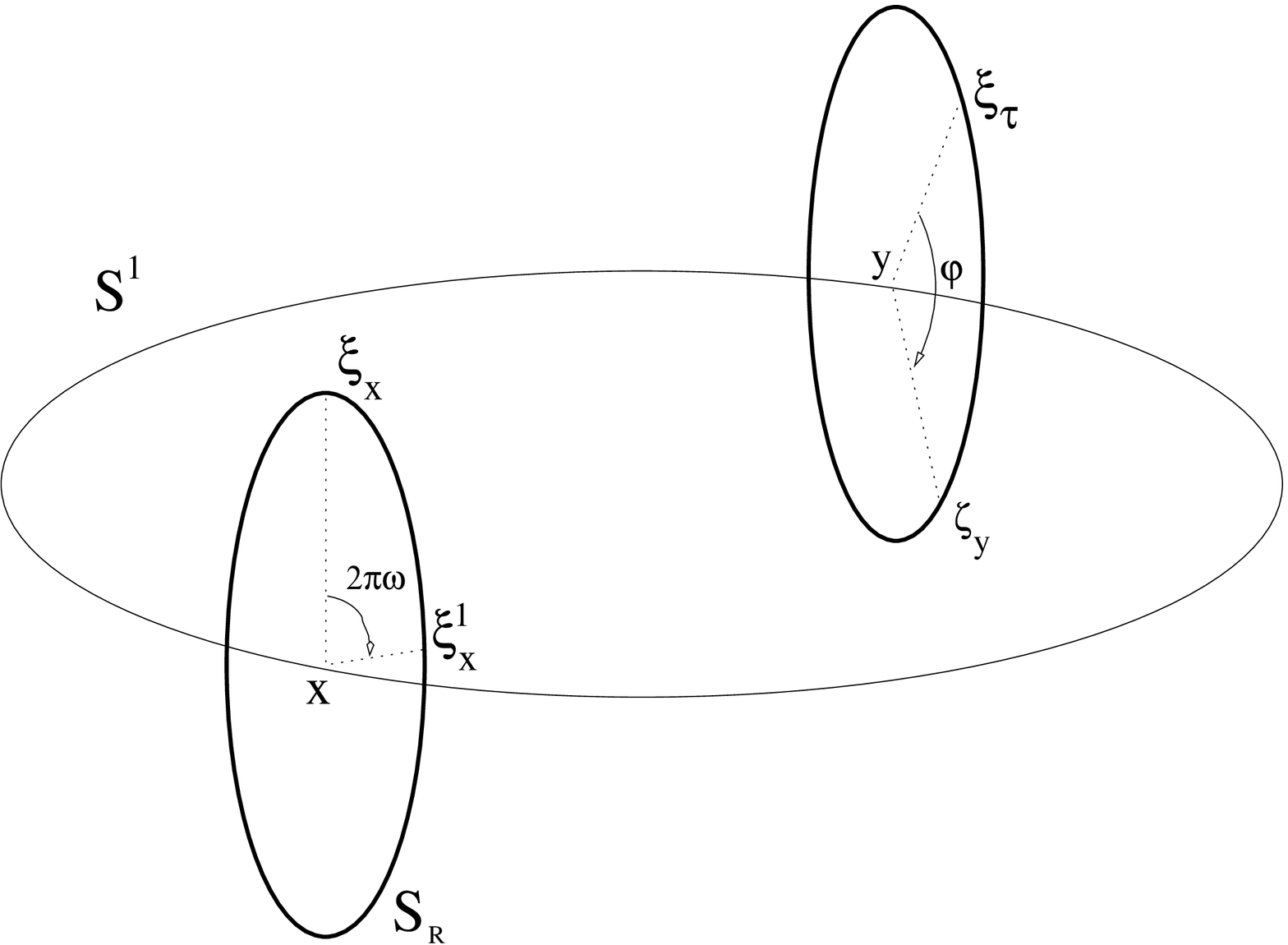}}}}
\end{center}

\subsection{Detail of the computation in the case $n=2$}

The following
subsection is technical, interpretation and discussion of the main result,
proposition \ref{propcarat} below,  are postponed to the next subsection.

To simplify the computation of the spectral distance it is useful
to delimit the part of $\aa$ that is relevant in the search of the supremum.
A first set of
simplifications comes 
from the low-dimension of the internal algebra.
{\lemma{\label{lemdiag1} The search for the supremum in computing
$d(\xox, \yoz)$ can be restricted to the set of selfadjoint
$a\in\aa$ whose diagonal $a_1$ and off-diagonal $a_0\doteq a- a_1$
are such that
\begin{eqnarray}
\label{premiereq}
&a_1(x) = 0, \;\tr{a_1(y)} \geq 0,&\\
&\label{posita1} \Delta_\xi(a_1) \geq 0,\quad \Delta_\xi(a_0)\geq
0&\\
&\label{aiizxi} z_\xi \delta(y) \geq 0&
\end{eqnarray}
where $z_\xi$ is defined in (\ref{defz}) and
\begin{equation}
\label{defdelta} \Delta_\xi\doteq \yoz - \xox,\quad \delta \doteq
a_{11} - a_{22}
\end{equation}
with $a_{ij}$ the components of $a$.}}
\newline

\noindent{\it Proof.} That $a$ is selfadjoint comes from a general
result according to which the supremum in the distance formula can
be searched on positive elements (see [\citelow{finite}], lemma 1). The
diagonal part $a_1$ enters the commutator condition,
\begin{equation}
\label{normdirect} \norm{[\dd,a]}= \underset{x\in S^1}{\sup}\lp
\frac{1}{2}\abs{\tr{\dot{a_1}(x)}} +
\sqrt{\frac{\dot{\delta}^2(x)}{4} + \norm{[\dd,a_0](x)}^2}\rp
\end{equation}
only via its derivative. Moreover for $\zeta_y\in\mathbb{T}_\xi$ and any constant $K_i$, $i=1,2$,
\begin{equation}
\label{defei}
\Delta_\xi (a_1 - K_i e_i) = \Delta_\xi (a_1)
\end{equation}
where $e_i$ is the projector on the $i^{th}$ term of the diagonal.
Therefore if $a$ attained the supremum, so does $a - K_i e_i$ whose
diagonal vanishes at $x$ as soon as one chooses $K_i = a_{ii}(x)$.
Hence the first term of (\ref{premiereq}).

Let $G$ be the group, acting on the components of $a$, generated by
the replacements
\begin{eqnarray}
\label{remplun}
a_1&\rightarrow& -a_1,\\
\label{rempldeux} a_0&\rightarrow& - a_0 \end{eqnarray} and the
permutation
\begin{equation}
\label{permtrois} a_{11}\leftrightarrow a_{22}. \end{equation}
 Since
(\ref{normdirect}) is invariant under the action of $G$ and
$a\rightarrow -a$ belongs to $G$, $\Delta_\xi(a)$ can be assumed
positive. Moreover $a$ attaining the supremum means that
\begin{equation}
\label{gupg}
\Delta_\xi(a) \geq \Delta_\xi(g(a)) \end{equation} for any $g\in G$.
For instance $g$ given by
(\ref{remplun}) yields
$$\Delta_\xi(a) = \Delta_\xi(a_1) +
\Delta_\xi(a_0 )  \geq -\Delta_\xi(a_1) +
\Delta_\xi(a_0 )$$ hence the first term of (\ref{posita1}). The second
term of (\ref{posita1}) is obtained putting  $g=$(\ref{rempldeux}) in
(\ref{gupg}). Similarly (\ref{aiizxi}) comes from
(\ref{permtrois}). $\tr{a_1(y)}$ being positive comes from
(\ref{remplun}) followed by (\ref{permtrois}), using the vanishing of $\xi_x(a_1)$. \hfill $\blacksquare$
\newpage
\noindent Other simplifications come from the choice of $S^1$ as the
base manifold.
{\lemma{
 \label{lemdiag2} Let $a= a_0 + a_1$ with components $a_{ij}$ as in lemma above. If
$a$ satisfies the commutator norm condition, then for $i=1,2$
\begin{equation}
\label{aiif}
 \abs{a_{ii}(\tau)} \leq F(\tau)
\end{equation}
where $F(\tau)$ is the $2\pi$-periodic function defined on $[0,
2\pi[$ by
\begin{equation}
\label{F}
 F(\tau)\doteq \text{min}\lp \tau, 2\pi - \tau\rp.
\end{equation}
Meanwhile
\begin{equation}
\label{azero} a_0=\dm{cc} 0& ge^{-i\Theta} \\
\overline{g}e^{i\Theta} & 0 \fm
\end{equation}
where $$\Theta\doteq \Theta_1 -\Theta_2$$ ($\Theta_i$ is defined in
(\ref{gdtheta})) and $g$ is a smooth function on $\rr$ given by
\begin{equation}
\label{g2} g(\tau) = g_0 + \int_{0}^\tau \rho(u)e^{i\phi(u)}du
\end{equation}
with $\rho\in C^\infty(\rr,\rr^+)$, $\norm{\rho}\leq 1$, and
$\phi\in C^\infty(\rr,\rr)$ satisfying
\begin{equation}
\label{krho} \rho(u+2\pi)e^{i\phi(u+2\pi)} =
\rho(u)e^{i(\phi(u)+\Theta(2\pi))}
\end{equation}
while the integration constant $g_0\in\cc$ is fixed by the equation
\begin{equation}
\label{gperiode}  g_0(e^{i\Theta(2\pi)} - 1)= \int_{0}^{2\pi}
\rho(u)e^{i\phi(u)}du.
\end{equation}
}}

\noindent{\it Proof.} (\ref{aiif}) comes from the commutator norm
condition $$\norm{\dot{a}_{ii}} = \norm{e_i [\dd,a]
e_i}\leq\norm{[\dd,a]}$$ where $e_i$ has been defined below (\ref{defei}), together with the
$2\pi$-periodicity of $a_{ii}$ (\ref{aij}) namely
$$a_{ii} = \int_{0}^\tau \dot{a}_{ii}(u) du = - \int_{\tau}^{2\pi}
\dot{a}_{ii}(u) du.$$

 The explicit form of $a_0$ is obtained by
noting that any complex smooth function $a_{12}\in\aa_E$ can be
written $ge^{-i\Theta}$ where $\Theta = \Theta_{12}$ is defined in
(\ref{thetaij}) and
 $g\doteq
a_{12}e^{i\Theta }\in C^\infty(\rr)$. By (\ref{lineairetheta})
\begin{equation}
\label{g} g(\tau + 2\pi) = g(\tau)e^{i\Theta(2\pi)}.
\end{equation}
Therefore any selfadjoint $a_0$ writes as in (\ref{azero}), which yields
by (\ref{diraac})
\begin{equation}
\label{commutazero} i[{\mathcal{D}},a_0] =\left(
\begin{array}{cc}
0 & \dot{g} e^{-i\Theta}\\ \dot{\overline{g}}e^{i\Theta}&
0\end{array}\right).
\end{equation}
By (\ref{normdirect}) the commutator norm condition implies
$$\norm{[{\mathcal{D}},a_0]}= \norm{\dot{g}}\leq 1,$$ that is to say
\begin{equation}
\label{g2bis}
 g(\tau)= g(0) + \int_{0}^\tau \rho(u)e^{i\phi(u)}du
\end{equation}
where $\rho\in C^{\infty}(\rr,\rr^+)$, $\norm{\rho}\leq 1$, $\phi\in
C^\infty(\rr,\rr)$. (\ref{gperiode}) is obtained from (\ref{g2})
inserted in (\ref{g}) for $\tau= 0$. Finally (\ref{gperiode})
inserted back in (\ref{g}) gives $$\int_0^{\tau}
\rho(u+2\pi)e^{i\phi(u+2\pi)} du = \int_0^\tau \rho(u)
e^{i(\phi(u)+\Theta(2\pi))} du,
$$ hence (\ref{krho}) by derivation with respect to $\tau$.\hfill $\blacksquare$
\\

\begin{figure}[h*]
\begin{center}
\mbox{\rotatebox{0}{\scalebox{.7}{\includegraphics{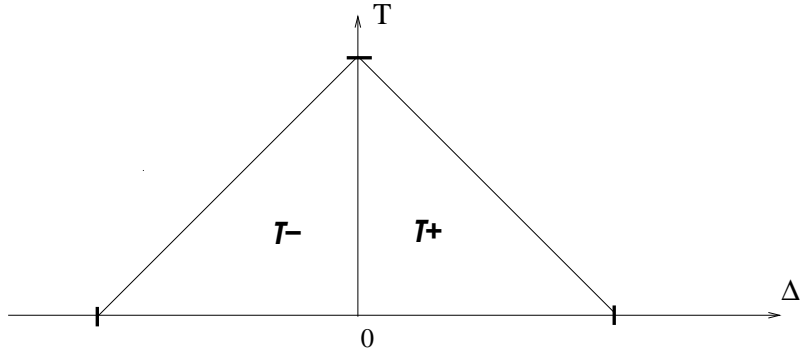}}}}
\end{center}
\caption{\label{tplusmoins} The diagonal part of $a$. Unit is
$\text{min}(\tau_0, 2\pi-\tau_0)$. }
\end{figure}

\noindent In [\citelow{cc}] we obtained by simple algebraic manipulations that 
for equatorial states the supremum could be searched on elements with vanishing diagonal.
Similarly for two states on the same fiber it is easy to see from  (\ref{premiereq}) that the diagonal 
part does not enter into play. For non equatorial states on different fibers however, the diagonal elements play a non trivial role,
determined by the value of $z_\xi$. 
{\corollary{\label{tplusmoinscorrol} The search for the supremum in
the computation of $d(\xox, \yoz)$ can be restricted to elements $a$
whose diagonal components at $y$,
\begin{equation}\label{tdelta}
\Delta \doteq \frac{1}{2}\delta(y),\quad
T \doteq \frac{1}{2}\tr a(y)\end{equation} lie within the triangle
\begin{equation}
\mathcal{T}_\pm \doteq T  \pm \Delta \leq \text{min}(\tau_0, 2\pi
-\tau_0)\end{equation} where the sign is the one of $z_\xi$ (fig.
\ref{tplusmoins}).}}
\\

\noindent{\it{Proof.}} When $z_\xi\geq 0$, $\Delta$ is positive by
(\ref{aiizxi}) (i.e. $a_{11}(y)\geq a_{22}(y)$) so that at least
$a_{11}(y)$ is positive otherwise (\ref{premiereq}) does not hold.
Therefore 
$$T + \Delta = a_{11}(y) \leq F(\tau_0)$$ by (\ref{aiif})
and the parametrization (\ref{parameters}) of $\yoz$. Similarly for
negative $z_\xi$ at least $a_{22}(y)$ is positive, hence $T - \Delta
= a_{22}(y) \leq F(\tau_0)$.\hfill$\blacksquare$


\newpage

 Thanks to these preliminary results it is
easy to come to the main result of this section, namely the
computation of all the distances on $\mathbb{T}_\xi$.

\begin{proposition}\label{propcarat} Let $\xox$ be a pure state in $P$ and $\yoz$
 a pure state in ${\mathbb{T}_\xi}$ parametrized according to (\ref{parameters})
 by a triple $(k, \tau_0, \varphi)$. Then either
 the two directions are far from each other
so that $\text{Con}(\xox)=\text{Acc}(\xox)$ and
\begin{equation}
\label{caratzero} \displaystyle d(\xox, \yoz)= \left\{
\begin{array}{ll}
\text{min} (\tau_0, 2\pi - \tau_0) &\text{when } \varphi = 0\\
+\infty &\text{when } \varphi \neq 0
\end{array}
\right.;
\end{equation}
or the directions are close to each other so that
$\text{Con}(\xox)={\mathbb{T}}_\xi$ and
\begin{equation}
\label{carat} \displaystyle d(\xox, \yoz)=
\underset{\mathcal{T}_\pm}{\max} \; H_\xi(T,\Delta)
\end{equation}
where the sign is the one of $z_\xi$ and
\begin{equation}
\label{defH} H_\xi(T,\Delta)\doteq T + z_\xi\Delta +
RW_{k+1}\sqrt{(\tau_0-T)^2-\Delta^2}
+RW_{k}\sqrt{(2\pi-\tau_0-T)^2-\Delta^2}
\end{equation}
with
\begin{equation}
\label{wmaxx} W_k  \doteq
\frac{\abs{\sin(k\omega\pi+\frac{\varphi}{2})}}{\abs{\sin
\omega\pi}}.
\end{equation}
\end{proposition}

\noindent{\it Proof.} (\ref{caratzero}) comes from a general result
[\citelow{cc}, Prop. 3] according to which for trivial holo\-nomy
the spectral distance coincides with the horizontal distance, which
in turn is either infinite or coincides with the geodesic distance
on the base. In the present context one can retrieve this result by
explicit calculation and we comment on this point in the paragraph
next to this proof. For the time being, we assume $\omega\neq 0$.
Let $a$ be an element of $\aa$ given by lemma \ref{lemdiag1}. By
(\ref{posita1}),
\begin{equation}
\label{redeux} \abs{\Delta_\xi(a)} = \Delta_\xi(a_0) +
\Delta_\xi(a_1)
\end{equation}
so that a good strategy to find the supremum is to begin with
constraining $\Delta_\xi(a_0)$ on the one side
(eq.(\ref{boundzeroxzero}) below), then adding $\Delta_\xi(a_1)$.
Writing
\begin{equation}
\label{defr} 2V_1\overline{V_2} = Re^{i\theta_0},
\end{equation}
Eqs. (\ref{evaluationdeux}) and (\ref{evaluationun}) written for
$\xi_x$ and $\yoz$ (the latest given by (\ref{zetay})) yield
$$\Delta_\xi(a_0) = \Re\lp{R e^{-i\theta_0}(g(\tau)e^{i\varphi} - g(0))}\rp,$$
that we rewrite, using (\ref{g2}),
\begin{equation}
\label{interm} \Delta_\xi(a_0) = R\int_0^\tau \rho(u) \cos \phi'(u)
+ \Re \lp R e^{-i\theta_0} g_0(e^{i\varphi} - 1)\rp,
\end{equation}
where
\begin{equation}
\label{dephfi} \phi'(u)\doteq \phi(u) -\theta_0 +\varphi.
\end{equation}
By the definition (\ref{tauk}) of $\tau$, using (\ref{krho}), the
integral term of (\ref{interm}) splits into
\begin{equation}
\label{splitunzero} \Re \int_0^{2k\pi} \rho(u) e^{i\phi'(u)}du = \Re
\lp\underset{n=0}{\overset{k-1}{\Sigma}} e^{2in\omega\pi}
\int_0^{2\pi} \rho(u) e^{i\phi'(u)}du\rp
\end{equation}
and
\begin{equation}
\label{splitdeuxzero} \Re \int_{2k\pi}^{2k\pi + \tau_0} \rho(u)
e^{i\phi'(u)}du = \Re \lp e^{2ik\omega\pi}\int_{0}^{\tau_0} \rho(u)
e^{i\phi'(u)}du\rp
\end{equation}
that recombine as
\begin{equation}
\label{splitun} S_{k+1}\int_{0}^{\tau_0} \rho(u) \cos \phi_k(u)\,du
+ S_{k}\int_{\tau_0}^{2\pi} \rho(u) \cos \phi_{k-1}(u)\,du
\end{equation}
where
\begin{equation}
\label{defsk} S_k\doteq \frac{\sin k\omega\pi}{\sin \omega\pi}, \,
\text{ and } \,\phi_k = \phi' + k\omega\pi.
\end{equation}
Meanwhile, after a few calculations with (\ref{gperiode}), the
real-part term of (\ref{interm}) writes
\begin{equation}
\label{splitdeux} S_{\frac{1}{2}}\int_0^{2\pi}\rho(u)\, \text{cos}
\, \phi_{\frac{1}{2}}(u)du
\end{equation}
where
\begin{equation} \label{sundemi} S_{\frac{1}{2}} \doteq
\frac{\sin \varphi/2}{\sin \,\omega\pi} \, \text{ and }\,
\phi_{1/2}\doteq \phi' - \frac{\varphi}{2} -\omega\pi.
\end{equation}
Combining (\ref{splitun}) and (\ref{splitdeux}) one finally obtains
\begin{equation}
\label{genquatre} \Delta_\xi(a_0)\leq R W_{k+1}
  \int_{0}^{\tau_0} \rho(u) \,du  +
  R W_{k}\int_{\tau_0}^{2\pi} \rho(u)\, du
\end{equation}
where $W_k$ is the maximum on $[0,2\pi[$ of $\abs{G_k(u)}$ defined
by
\begin{equation}
\label{gk} G_k\doteq S_k\cos \phi_{k-1} +
S_\frac{1}{2}\cos\phi_\frac{1}{2}.
\end{equation}
So far we used the commutator norm condition on $a_0$ and $a_1$
independently, via lemmas \ref{lemdiag1} and \ref{lemdiag2}. To
obtain a global constraint on $a$, first note that
\begin{equation}
\norm{[\dd,a_0]} = \norm{-i\dm{cc}0 & \dot{g}e^{-i\Theta}\\
\dot{\overline{g}}e^{i\Theta} & 0 \fm} = \rho
\end{equation}
by (\ref{g2}), so that (\ref{normdirect}) yields for any $a$
satisfying the commutator norm condition
\begin{equation} \label{normdelata} \rho \leq \sqrt{\lp
1-\frac{\abs{\text{Tr}\,\dot{a_1}}}{2}\rp^2 -
\lp\frac{\dot{\delta}}{2}\rp^2}
\end{equation}
that is, by Jensen inequality,{\footnote{For $a\leq b$ and any $f\in
L_1([a,b])$ ($f$ positive in the first equation) one has\cite{rudin}
 \begin{eqnarray}
 \label{jensenun}
\int_a^b \sqrt{f(t)}dt&\leq& \sqrt{b-a}\,\sqrt{\int_a^b f(t)dt},\\
\int_a^b f^2(t)dt&\geq& \frac{1}{b-a}\lp\int_a^b f(t) dt\rp^2,\\
\int_a^b \lp 1-f(t)\rp^2dt &\leq& (b-a)\lp 1 -\frac{1}{b-a}\int_a^b
f(t) dt\rp^2.
\end{eqnarray}
}}
\begin{equation}
\label{mink1}
 \int_a^b\rho(t)dt\leq \sqrt{\lp (b-a) - \frac{\abs{\tr{a_1}(a) - \tr{a_1}(b)}}{2} \rp^2 -
 \lp\frac{\delta(a) -
 \delta(b)}{2}\rp^2}.
\end{equation}
Therefore (\ref{genquatre}) together with the vanishing of $a_1$ at
$x$ and the positivity of $\tr{a_1(y)}$ yields the announced bound
\begin{equation}
\label{boundzeroxzero} \Delta_\xi(a_0)\leq RW_{k+1}\sqrt{\lp\tau_0 -
T\rp^2 - \Delta^2} + RW_{k}\sqrt{\lp (2\pi -\tau_0) - T\rp^2-
\Delta^2}
\end{equation}
where we use the notations (\ref{tdelta}). Adding
\begin{equation}
\label{deltaxia10}
 \Delta_\xi(a_1)
  = T +  z_\xi\Delta
\end{equation}
one obtains that $\Delta_\xi(a)$ is inferior or equal to the right
hand side of (\ref{carat}).

To prove the r.h.s. of (\ref{carat}) is the lowest upper bound, we
need a sequence of elements
\begin{equation}
\label{an} a_n = \dm{cc} f_n^+ & g_n e^{-i\Theta}\\ \overline{g}_n
e^{i\Theta}&f_n^-\fm
\end{equation} such that
\begin{equation}
\label{rhscarat} \underset{n\rightarrow +\infty}{\text{lim}}
\Delta_\xi(a_n) = \underset{\mathcal{T}_\pm}{\text{max}} \;
H_\xi(T,\Delta).
\end{equation}
The diagonal part of $a_n$ is defined by the values
$(T_0,\Delta_0)\in\mathcal{T}_\pm$ for which $H$ attains its maximum
(see corollary \ref{corollcarat} below) in the following way:
assuming $\tau_0\neq 0$, then $f_n^{\pm}$ is a sequence of smooth
functions approximating from below the $2\pi$-periodic continuous
function,
\begin{equation}
\label{fmoinsplus} \;f^\pm(t) = \left\{
\begin{array}{ll}
C^\pm t&\text{ for } 0\leq t\leq \tau_0\\
C^\pm \tau_0 - \frac{C^\pm\tau_0}{2\pi-\tau_0}(t-\tau_0)&\text{ for
} \tau_0\leq t\leq 2\pi
\end{array}\right. \text{ with } C^\pm\doteq \frac{T_0 \pm
\Delta_0}{\tau_0}.\end{equation} In other words we build $a_{n}$ in
such a way that
\begin{equation}
\label{trasol}
\underset{n\rightarrow+\infty}{{\lim}}\;\frac{1}{2}\tr a_{n}(t) =
\left\{
\begin{array}{ll}
\frac{T_0}{\tau_0} t&\text{ for } 0\leq t\leq \tau_0\\
T_0  - \frac{T_0}{2\pi-\tau_0}(t-\tau_0)&\text{ for } \tau_0\leq
t\leq 2\pi
\end{array}\right.\end{equation}
while, writing $\delta_n\doteq f_n^+ - f_n^-$,
\begin{equation}
\label{deltasol}
\underset{n\rightarrow+\infty}{{\lim}}\;\frac{1}{2}\delta_n(t) =
\left\{
\begin{array}{ll}
\frac{\Delta_0}{\tau_0} t&\text{ for } 0\leq t\leq \tau_0\\
\Delta_0 - \frac{\Delta_0}{2\pi-\tau_0}(t-\tau_0)&\text{ for }
\tau_0\leq t\leq 2\pi
\end{array}\right. .\end{equation}
Thus the diagonal part $a_{1,n}$ of $a_n$ does satisfy lemmas
\ref{lemdiag1},\ref{lemdiag2} and yields
\begin{equation}
\label{evaltracesol} \underset{n\rightarrow+\infty}{{\text lim}}
\Delta_\xi(a_{n,1}) = T_0 + z_\xi\Delta_0.\end{equation} In case of
vanishing $\tau_0=0$, one simply choose $a_{n,1}=0$. The
off-diagonal part $a_{n,0}$ of $a_n$ is defined by substituting
$\rho$ and $\phi$ in the definition (\ref{g2bis}) with sequences
$\rho_n$, $\phi_n$ approximating from below the $2\pi$-periodic step
functions $\Gamma$, $\Phi$ defined as follows: $\Gamma$ is given by
(\ref{normdelata}) with $\tr a_1$ and $\delta$ replaced by their
limit values above, namely
\begin{equation}
\Gamma\doteq \left\{\begin{array}{ll} \sqrt{ \lp
1-\frac{T_0}{\tau_0}\rp^2 -
\lp\frac{\Delta_0}{\tau_0}\rp^2} &\text{ for } 0\leq t\leq \tau_0\\
 \sqrt{ \lp
1-\frac{T_0}{2\pi - \tau_0}\rp^2 - \lp
\frac{\Delta_0}{2\pi-\tau_0}\rp^2}&\text{ for } \tau_0\leq t\leq
2\pi
\end{array}\right. .\end{equation}
In case $\tau_0=0$, one takes $\Gamma = 1$. This guarantees that
$a_n$ satisfy the commutator norm condition. In order to define
$\Phi$, one needs to work out $W_k$ explicitly. Easy calculations
from $(\ref{gk})$ yields
$$G_k = A_k \cos \phi' + B_k\sin \phi' \text{ with } \left\{
\begin{array}{l}
A_k\doteq S_\frac{1}{2}\cos(\frac{\varphi}{2}+\omega\pi) + S_k\cos(k-1)\omega\pi,\\
B_k\doteq S_\frac{1}{2}\sin(\frac{\varphi}{2}+\omega\pi) -
S_k\sin(k-1)\omega\pi,
\end{array}\right.$$
so that $G_k$ attains its maximum value
\begin{equation}
\label{wmax} W_k  \doteq \sqrt{A_k^2 + B_k^2} =
\frac{\abs{\sin(k\omega\pi+\frac{\varphi}{2})}}{\abs{\sin(\omega\pi)}}
\end{equation}
when $$\text{tan}\, \phi' = \frac{B_k}{A_k}$$ that is to say, after a
few algebraic manipulations, when
$$\cos^2 \phi' =
\frac{A_k^2}{W_k^2}= \cos^2\lp(k-1)\omega\pi-\frac{\varphi}{2}\rp.$$
Thus, remembering (\ref{dephfi}), $G_k[u]$ is constantly maximum on
$[0,2\pi]$ as soon as
{\footnote{note the change of notation with respect to [\citelow{cc}]
where $\Phi_k$ was denoting $\text{Argtan}\frac{A_k}{B_k}$}}  $$\phi(u) = \Phi_k$$ where
\begin{equation}
\label{phimaxbis} \Phi_k \doteq (1-k)\omega\pi -\frac{\varphi}{2} +
\theta_0\, \text{ or } \, \Phi_k \doteq (1-k)\omega\pi
-\frac{\varphi}{2} + \theta_0 + \pi
\end{equation}
depending on the sign of $A_k$ and $B_k$.
Hence the required $2\pi$-periodic step function $\Phi$ defined on
$[0,2\pi]$ by
\begin{equation}
\label{phi} \Phi(u) =\left\{ \begin{array}{ll} \Phi_{k+1}
&\text{ for } 0 \leq u < \tau_0\\
\Phi_{k} &\text{ for } \tau_0 < u < 2\pi.
\end{array} \right.
\end{equation}
Repeating the procedure that leads to (\ref{genquatre}), on easily
checks that
\begin{eqnarray}\underset{n\rightarrow+\infty}{\lim}\Delta_\xi(a_{n,0})&=& RW_{k+1}\int_0^{\tau_0} \Gamma(u) du +
RW_{k}\int_{\tau_0}^{2\pi} \Gamma(u) du\\
&=& RW_{k+1}\sqrt{ \lp \tau_0-T_0\rp^2 - \Delta_0^2} + RW_{k}\sqrt{
\lp 2\pi-\tau_0 -T_0\rp^2 -  \Delta_0^2}\end{eqnarray} Adding
(\ref{evaltracesol}), one finally gets
$\underset{n\rightarrow+\infty}{\
lim}\Delta_\xi(a,n)=H_\pm(T_0,\Delta_0)$ as expected in
(\ref{rhscarat}). Hence (\ref{carat}). \hfill $\blacksquare$
\\

One could be tempted to search a more explicit formula of
$d(\xox,\yoz)$ by determining the maximum of the function $H$. This
is not clear whether this maximum has such a simpler form. To see
where is the difficulty, let us first examine the easy case of
vanishing $\omega$. The integral part of (\ref{interm}) is still
given by (\ref{splitun}) assuming $S_k=k$. Regarding the real part
term of (\ref{interm}), either it can be made infinite when
$\varphi$ does not vanish (since (\ref{gperiode}) does not constrain
$g_0$ any longer) thus yielding the second line of
(\ref{caratzero}), or it vanishes when $\varphi = 0$ and
(\ref{splitdeux}) makes sense assuming $S_{\frac{1}{2}} = 0$.
Practically one obtains
\begin{equation}
\label{boundzeroxzerox} \Delta_\xi(a_0)\leq R(k+1)\int_{0}^{\tau_0}
\rho(u) \cos \phi'(u)\,du + Rk\int_{\tau_0}^{2\pi} \rho(u) \cos
\phi'(u)\,du.
\end{equation}
Now the vanishing of $\omega\pi=0$ makes the function $g$ defined in
(\ref{g}) $2\pi$-periodic, therefore
\begin{equation}
\int_0^{2\pi} \rho(u)e^{i\phi(u)} =0
\end{equation}
and (\ref{boundzeroxzerox}) writes
\begin{eqnarray}
\label{boundzeroxzeroxx} \Delta_\xi(a_0)&\leq& R\int_{0}^{\tau_0}
\rho(u) \cos \phi'(u)\,du = -R\int_{\tau_0}^{2\pi} \rho(u) \cos
\phi'(u)\,du\\
&\leq& R\,\text{min} \lp \int_{0}^{\tau_0} \rho(u) \,du ,
\int_{\tau_0}^{2\pi} \rho(u) \,du\rp.
\end{eqnarray}
By Jensen inequalities one has, instead of the r.h.s. of
(\ref{carat}),
\begin{equation}
\label{boundzeroxzeroproof} \Delta_\xi(a)\leq \text{min}\lp
\underset{\mathcal{T}_\pm}{\max} \,H_{\tau_0},\;
\underset{\mathcal{T}_\pm}{\max}\, H_{2\pi-\tau_0}\rp
\end{equation}
where
\begin{equation}\label{hlambda}
H_\lambda \doteq T + z_\xi\Delta + R\sqrt{(\lambda- T)^2 -
\Delta^2}.
\end{equation}
The maximum of $H_\lambda$ is easily found by noting that for a
fixed value of $T$, $H_\lambda$ is maximum for $\Delta^2 =
z^2_\xi(\delta - T)^2$, which yields $H_\lambda = \lambda$ and
(\ref{hlambda}) = (\ref{caratzero}) as expected.

For non vanishing $\omega$ the situation is more complicated since
$H_\xi$, unlike $H_\lambda$, involve two square roots. In fact there
does not seem to be a simple form for the maximum of $H_\xi$. So far
the best we managed to do is summarized in the following corollary.

{\corollary{\label{corollcarat}$H_\xi$ reaches its maximum either on
the segment $T=0$ or on the segment $T\pm \Delta = \min(\tau_0,2\pi
-\tau_0$).}}
\\

\noindent{\it Proof.} Let us first show that $H_\xi$ attains it
maximum on the border of $\mathcal{T}_\pm$. A local maximum
$p=(T,\Delta)$ of $H_\xi$ inside ${\mathcal T}_\pm$ is a solution of
the system
\begin{equation}
\label{systemeun} {\frac{\partial H_\xi}{\partial T}}_{\slash p} =
0, \; \frac{\partial H_\xi}{\partial \Delta}_{\slash p} = 0.
\end{equation}
Seen as a system in the non-vanishing variables
\begin{equation}
\label{defivi} V_1 =\sqrt{(\tau_0-T)^2 -\Delta^2},\, V_2
=\sqrt{(2\pi -\tau_0-T)^2 -\Delta^2},
\end{equation}
 (\ref{systemeun} has solution
\begin{equation} \label{systeme} \left\{ \begin{array}{ccl}
\Delta - Z(\tau_0-T) &=& 2RW_k\Delta V_2^{-1}(\pi
- \tau_0),\\
Z(2\pi - \tau_0-T) -\Delta &=& 2RW_{k+1}\Delta V_1^{-1}(\pi -
\tau_0).
\end{array}\right. \end{equation}
For $\pi = \tau_0$, the above system in $T,\Delta$ is solved by the
straight line $l$
\begin{eqnarray}
\label{solpizzero} \Delta = 0 &\text{ when } Z = 0,\\
\label{solpiznonzero}T = \pi - \frac{\Delta}{Z} &\text{ when } Z
\neq 0.
\end{eqnarray}
The latest yields
\begin{equation}
H_\xi(T,\Delta)= \pi + \frac{\Delta R^2}{Z}(W_k + W_{k+1} -1)
\end{equation}
whose derivative with respect to $\Delta$ is constant. Therefore
whatever $Z$, $H_\xi$  attains its maximum on one of the end point
of $l\cap\mathcal{T}_\pm$, i.e. $H_\xi$ is maximum on the border of
$\mathcal{T}_\pm$. For $\pi \neq \tau_0$, equaling $V_1$ and $V_2$
expressed by (\ref{systeme}) to their definition (\ref{defivi}), one
finds that $p$ is a common root of two polynomials
\begin{eqnarray*}
P_1 &=& \Delta^4-2Z\Delta^3(2\pi-\tau_0-T)+\Delta ^2(4R^2W_{k+1}^2
(\pi-\tau_0)^2+Z^2
(2\pi-\tau_0-T)^2-(\tau_0-T)^2) \\
&&+2Z\Delta(2\pi-\tau_0 -T)(\tau_0-T)^2-Z^2
(2\pi-\tau_0-T)^2(\tau_0-T)^2,\\
P_2 &=&\Delta^4-2Z\Delta^3(\tau_0-T)+\Delta^2(4R^2W_k^2(\pi
-\tau_0)^2-(2\pi-\tau_0-T)^2+Z^2(\tau_0-T)^2)
\\
&&+2Z\Delta(2 \pi-\tau_0-T)^2(\tau_0-T) -Z^2(2 \pi-\tau_0-T)^2 (\tau
_0-T)^2.
\end{eqnarray*}
Since two polynomials have common roots if and only if their
discriminant vanishes, the coordinate $\Delta$ of $p$ is a root of
the resultant of $P_1$, $P_2$ viewed as polynomials in $T$. Thanks
to formal computation programs{\footnote{Here we used Mathematica}}
one finds that this resultant is the product of 
$$256 Z^8\Delta^8(\pi
- \tau_0)^6$$ by
\begin{eqnarray*}
P_3\doteq R^4\Delta^2 &-& 4R^2Z\Delta(\tau_0-\pi)
\frac{\sin((2k+1)\omega\pi+\varphi)}{\sin \omega\pi}\\
&+&4(\tau_0-\pi)^2\frac{(R^2\cos(k\omega\pi +
\varphi)^2-1)(R^2\cos((k+1)\omega\pi + \varphi)^2+1)}{\sin^2
\omega\pi}.
\end{eqnarray*}
The discriminant of $P_3$ viewed as an equation in $\Delta$ of the
second degree is $$-\frac{16R^4(\pi-\tau_0)^2}{\sin^2 \omega\pi}
\left( Z^2\cos(k\omega\pi+\varphi)\cos((k+1)\omega\pi+\varphi)
-\sin(k\omega\pi +\varphi)\sin((k+1)\omega\pi+\varphi) \right)^2,$$
hence $P_3$ has no real solution. In other terms, the function
$H_\xi$ does not have local extrema inside $\mathcal{T}_\pm$, hence
it reaches its maximum on the border of $\mathcal{T}_\pm$.

Now observing that
\begin{equation}
\label{defder}
 \frac{\partial H_\xi(T,0)}{\partial T}= 1 -RW_{k+1}
-RW_{k}
\end{equation}
is a constant, one deduces that $H_\xi(T,0)$ is maximum at one of
the end points of the segment $\Delta=0$.\hfill$\blacksquare$
\\

\noindent Let us underline that the derivative of $H_\xi$ is
constant neither on the segment $T=0$ nor on the hypotenuse
$T\pm\Delta=\min(\tau_0,2\pi-\tau_0)$. For instance one can show
that $H_\xi(0,\Delta)$ has a maximum $h_1$ for a value $\Delta_0$
which is a root of a polynomial of degree $4$. We did not manage to
find a simpler expression for $\Delta_0$, neither we did for the
maximum $h_2$ of $H_\xi$ on the hypotenuse. Hence, unless
numerically specifying $k,\tau_0, z_\xi$ and $\varphi$, it seems
difficult to compare $h_2$ and $h_1$ and go beyond corollary
\ref{corollcarat}.
\newpage

\subsection{Discussion}

\subsubsection*{Equatorial versus non equatorial states}

In case $\xi$ is an equatorial state, i.e. $z_{\xi_{\text{eq}}} =
0$, corollary  \ref{corollcarat} and \ref{tplusmoinscorrol} indicates that $H_\xi$ attains
its maximum at the point $(T, \Delta)=(0,0)$. Thus
\begin{equation}
\label{revisun}
d(\xox,\yoz) =
X\doteq H_\xi(0,0)=  RW_{k+1}\tau_0 + RW_{k}(2\pi-\tau_0).
\end{equation}
To make the link with [\citelow{cc}], note that corollary \ref{tplusmoinscorrol}
alone yields
\begin{equation}
\label{revisedeux}
d(\xox,\yoz) =
\underset{\underset{\Delta=0}{(T,\Delta)\in\mathcal{T}_\pm}}{\max}
\left(X +TY\right)
\end{equation}
where
$$
 Y
\doteq \frac{\partial H_\xi(T,0)}{\partial T}$$
 given by
(\ref{defder}). This is nothing that proposition V.4
of [\citelow{cc}]. In the statement of this  proposition we
did not insist on the restriction to equatorial states{\footnote{In
the published version this point is however clearly indicated at the
beginning of the section, (eq. 61, op.cite). On the arXiv version
the restriction is recalled in the statement of the proposition.}}
(specifically we did not replace $R$ by its value $1$) because we
hoped the result would still be true for non vanishing $z_\xi$. Instead
 the present work shows that the spectral distance
depends on $z_\xi$ on a non-trivial
 manner. Now, putting $R = 1$ yields
$Y\leq 0$ since $\abs{\sin \omega\pi} \leq \abs{\sin (k+1)\omega\pi}
+\abs{\sin k\omega\pi}$ hence (\ref{revisun}). 

 For non equatorial states $Y$ can be positive and there is no
further conclusion than corollary \ref{corollcarat}.

\subsubsection*{The shape of the fiber}

 When $\tau_0=0$, i.e. when $\yoz=\zeta_x$ belongs to the
fiber of $\xox$, $\mathcal{T}_\pm$ reduces to the single point $0$
and
\begin{equation}
\label{dfibre}
 d(\xox,\yox) =  2\pi  RW_{k} = \frac{2\pi R}{\abs{\sin
\omega\pi}}\abs{\sin \frac{\Xi}{2}}
\end{equation}
where we parametrize the fiber over $x$ by
$$\Xi\doteq 2k\omega\pi + \varphi.$$
We retrieve formula (134) of [\citelow{cc}] which, as already mentioned
there, is valid also for non equatorial state. In [\citelow{smoother,african}]
we give an interpretation of this result saying that for an
intrinsic metric point of view, the $S^1$ fiber over $x$ equipped
with the spectral distance has the shape of a cardioid. In fact
(\ref{dfibre}) can also be viewed as the length of the straight
segment inside the circle{\setcounter{footnote}{0}\footnote{This has
been pointed out by the audience, during a talk given at Toulouse
university.}}. 
\begin{figure}[htbp]
\begin{center}
\mbox{\rotatebox{0}{\scalebox{.6}{\includegraphics{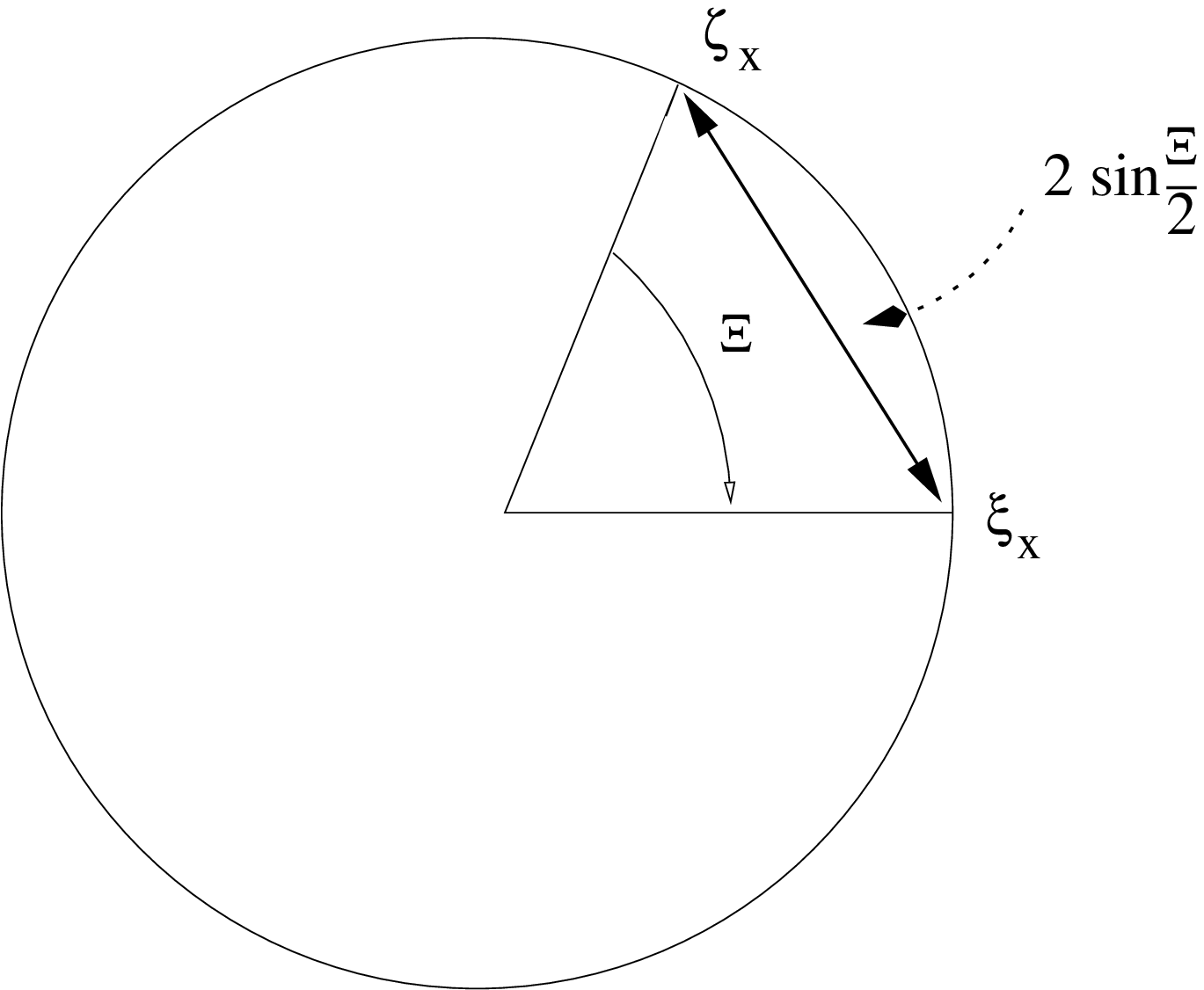}}}}
\label{figcercle}
\end{center}
\end{figure}
Therefore the spectral distance appears as the
geodesic distance {\it inside} the disk, in the same way that in the
two-sheets model\cite{cc} the spectral distance coincides with a geodesic
distance {\it within} the two sheets. Say differently, the spectral
distance coincides with the geodesic distance of a manifold $M'$
whose boundary is the set of pure states of $\aa$. Of course it is
appealing to identify $M'$ to the set of all states of $\aa$ (i.e.
all convex combinations of pure states) but this idea does not
survive here since the convex combination $\rho$ of $\xox,\yox$,
$$\left( \begin{array}{c}\rho_1\\ \rho_2\end{array}\right)
= \lambda \left( \begin{array}{c}\xi_1\\ \xi_2\end{array}\right) +
(1-\lambda) \left( \begin{array}{c}\zeta_1\\
\zeta_2\end{array}\right)$$ has an altitude $z_\rho$ that does not
equal $z_\xi= z_\zeta$. In other terms the straight segment $[\xox,
\xoz]$ inside $T_\xi$ does not coincides with the convex combination
$\lambda\xi + (1-\lambda\zeta)$. We deserve to further work a better
comprehension of this observation.

\section{Spectral distance on the circle with a $\cc P^{n-1}$ fiber}

Let us now investigate the general case
$$\aa=C^{\infty}(M,\mn)$$
for arbitrary integer $n\in\nn$. $\mathbb{T}_\xi$ is now an $n$-torus
and instead of (\ref{parameters}) one deals with equivalence classes
of $(n+1)$-tuples
\begin{equation}
\label{nparameters} \lp k\in \nn,\, 0\leq \tau_0\leq 2\pi, \, 0\leq
\varphi_i \leq 2\pi\rp \sim \lp k+ \zz,\, \tau_0,\, \varphi_j -
2\zz\omega_j\pi\rp
\end{equation}
with
\begin{equation}
\label{ntauk} \omega_j\doteq \frac{\Theta_1(2\pi) -
\Theta_j(2\pi)}{2\pi}\quad\quad\quad \forall j\in[2,n]
\end{equation}
such that $\yoz$ in $\mathbb{T}_\xi$ writes
\begin{equation}
\label{etattaun}
\yoz = \left(\begin{array}{r} V_1(\tau)\\
e^{i\varphi_j} V_j(\tau)\end{array}\right) \end{equation} where
$\tau\doteq 2k\pi + \tau_0$. As soon as $n>2$ there is no longer
correspondence between the fiber of $P$ and a sphere, however in
analogy with (\ref{defr}) we write
\begin{equation} \label{ndeltaxi}
V_j = \sqrt{\frac{R_j}{2}}e^{i\theta^0_j}
\end{equation}
where $R_j\in\rr^+$, $\theta^0_j\in[0,2\pi]$.

The simplifications relying on the choice of $S^1$ as a base,
namely lemma \ref{lemdiag2}, still hold for $n\geq 2$:

{\lemma\label{nlemmadiag}{Whatever $\yoz\in\mathbb{T}_\xi$, the
supremum in computing $d(\xox,\yoz)$ can be searched on selfadjoint
elements $a$ whose diagonal part $a_1$ vanishes as $x$ and whose
components $a_{ii}$ satisfy (\ref{aiif}). The non-diagonal
components $a_{ij}$ can be written
\begin{equation}
\label{offdiagn} a_{ij} = g_{ij}e^{-i\Theta_{ij}} \end{equation}
with $\Theta_{ij}$ defined in (\ref{thetaij}) and
\begin{equation}
\label{ndefgij} g_{ij}(\tau) = g_{ij}^0 + \int_0^\tau
\rho_{ij}(u)e^{i\phi_{ij}(u)}du \end{equation} where - for each
couple of indices $(i,j)$-  the real positive normalized smooth
function $\rho_{ij}$, the real smooth function  $\phi_{ij}$ and the
complex constant $g^0_{ij}$ satisfy equations (\ref{krho}) and
(\ref{gperiode}).}}
\\

\noindent{\it{Proof. }} The vanishing of $a_1$ and the boundary
values of the $a_{ii}$'s are obtained like in lemmas \ref{lemdiag1}
and \ref{lemdiag2}. The properties of the non-diagonal components
are obtained like in \ref{lemdiag2}, except that one has to begin
with $(e_i + e_j) a (e_i +e_j)$ instead of $a$, where $e_i$ denotes
the $i^{th}$ eigenprojector of the connection $1$-form $A$.\hfill
$\blacksquare$
\\

 For explicit computation, we
restrict in the following to a given fiber, say the one over $x\in
S^1$. In other words we consider only those $\yoz = \xoz$
given by (\ref{nparameters}) with $\tau_0=0$.

{\proposition  Given a pure state $\xoz = (k, 0,
\varphi_j)\in\mathbb{T}_\xi$, either $\xoz$ does not belongs to the
connected component $\mathbb{U}_\xi$ and
\begin{equation}
\label{infn} d(\xox,\xoz) = +\infty \end{equation} or
$\xoz\in\mathbb{U}_\xi$ and
\begin{equation}
\label{resultatnproof} d(\xox,\xoz) = \pi \tr \abs{S_k}
\end{equation}
where $\abs{S_k} = \sqrt{S_k^*S_k}$ and $S_k$ is the matrix with
components
\begin{equation}
\label{defSij} S^k_{ij}\doteq \sqrt{R_i R_j}\;\frac{\sin\lp
k\pi(\omega_j - \omega_i)+ \frac{\phi_j - \phi_i}{2}\rp}{ \sin
\pi(\omega_j - \omega_i)}
\end{equation}
where $\omega_i$ is defined in (\ref{ntauk}).}
\\

\noindent{\it Proof.} Staying within the same fiber, by lemma
\ref{nlemmadiag} the supremum can be searched on elements $a$ whose
diagonal part $a_1$ vanishes. Therefore (\ref{defdelta}) with $\yoz$ given by (\ref{etattaun}), together with (\ref{ndeltaxi})
and (\ref{ntauk}) yield
\begin{equation}
\label{ndeltaximoinsun} \Delta_\xi(a) = \sum_{i,j=1}^{n}
\frac{\sqrt{R_i R_j}}{2}\lp a_{ij}(x)
e^{i(\theta_{j}^0-\theta_{i}^0)}\lp e^{2ik\pi(\omega_j -\omega_i)}
e^{i(\varphi_j-\varphi_i)}-1\rp\rp.
\end{equation}
For any couple of directions $i,j$ far from each other,
$\omega_i=\omega_j$ so that the corresponding term in
(\ref{ndeltaximoinsun}) either vanishes if $\varphi_i=\varphi_j$, or
can be made infinite if $\varphi_i\neq\varphi_j$ since $g^0_{ij}$ is
no longer constrained by (\ref{gperiode}). In this case
$\xoz\notin\mathbb{U}_\xi$ by definition (see (\ref{txif})), hence
(\ref{infn}). Consequently (\ref{ndeltaximoinsun}) reduces to a sum
on couples of close directions. Then (\ref{ndefgij}) together with
(\ref{gperiode}) yield
\begin{equation}
a_{ij}(x) = g_{ij}^0 = \frac{1}{2i}\frac{e^{-i\pi(\omega_j
-\omega_i)}}{\sin \pi(\omega_j-\omega_i)}\int_0^{2\pi} \rho_{ij}(u)
e^{i\phi_{ij}(u)}du
\end{equation}
and (\ref{ndeltaximoinsun}) rewrites
\begin{equation}
\Delta_\xi(a) = \frac{1}{2}\sum_{i,j=1}^{n} S_{ij}^k\int_0^{2\pi}
\rho_{ij}(u) e^{i\phi^{k-1}_{ij}(u)} du
\end{equation}
where $S_{ij}^k$ is defined in (\ref{defSij}) and
\begin{equation}
\label{phiprimeij} \phi^k_{ij}\doteq \phi_{ij} + k\pi(\omega_j -
\omega_i) + \theta_{j}^0-\theta_{i}^0+\frac{\varphi_j -
\varphi_i}{2}.
\end{equation}
Calling $\Gamma$ the matrix with components
$\rho_{ij}e^{i\phi^{k-1}_{ij}}$,
\begin{equation}
\label{tracenn}
\Delta_\xi(a) =\frac{1}{2}\int_0^{2\pi} \tr(S^k
\Gamma(u)).
\end{equation}
For $a$ whose vanishing diagonal and off-diagonal given by
(\ref{offdiagn}) the commutator norm condition is easily computed
thanks to (\ref{diraac}),
\begin{equation}
\label{commutnexact} \norm{[\dd,a]}=\norm{\lp\begin{array}{cc} 0&
\rho_{ij}e^{i\phi_{ij}}e^{i(\Theta_j-\Theta_i)}\\
 \rho_{ij}e^{-i\phi_{ij}}e^{-i(\Theta_j-\Theta_i)}&0
 \end{array}\rp}\leq 1.
 \end{equation}
Introducing the diagonal matrix $E$ with components
$$E_{ii} = e^{i(\Theta_i - (k-1)\pi\omega_i - \theta_0^i - \frac{\varphi_i}{2})},$$
one observes that
$$\Gamma =  E [\dd,a] E^*$$
hence $\norm{\Gamma(u)}\leq 1$ for any $u\in[0,2\pi]$. Back to
(\ref{tracenn}),{\footnote{Writing $\alpha,\beta$ the components of
the matrices within the basis that diagonalizes $\Gamma$,
$$\tr\lp S_k\Gamma(u)\rp = \sum_{\alpha} S_{\alpha\alpha}^k
\Gamma_{\alpha\alpha}(u) \leq \sum_{\alpha}
\abs{S_{\alpha\alpha}^k}\norm{\Gamma}= \tr\abs{S_k}
$$}}
$$\tr\lp S_k\Gamma(u)\rp \leq \tr\abs{S_k}.$$
Hence
\begin{equation}
\label{resultatn} \Delta_\xi(a) \leq \pi\tr\abs{S_k}.
\end{equation}
This upper bound is reached by the diagonal constant matrix $\Gamma$
whose components in the basis that diagonalizes $S_k$ (with
corresponding indices $\mu,\nu$) are
$$
\Gamma_{\mu\mu} = \text{sign}\,S^k_{\mu\mu}
$$
and $\phi_{ij}$ the constant function
\begin{equation}
\label{phiijtest}\phi_{ij}(u) = (k-1)\pi(\omega_i-\omega_j) +
\theta_{i}^0 - \theta_{j}^0 + \frac{\varphi_i - \varphi_j}{2}.
\end{equation}
\hfill $\blacksquare$

To check consistency with the $n=2$ case, note the following
correspondence
\begin{eqnarray}
\omega \text{ in
(\ref{tauk})} &\text{reads}&  \omega_2 \text{ in (\ref{ntauk}) },\\
\varphi\text{ in (\ref{parameters})}&\text{reads}& \varphi_2 \text{
in (\ref{etattaun})},\\ \theta_0\text{ in (\ref{defder}) }
&\text{reads}& \theta^0_1 - \theta_0^2\text{ in
(\ref{ndeltaxi})},\\
\omega_1 = \varphi_1 = 0&\text{and}&R_1 = R_2 = R,
\end{eqnarray}
the latest coming from the requirement that the directions in the
$n=2$ case are close to each other. Within this correspondence, all
the quantities without indices used in the proof of of proposition
\ref{propcarat} corresponds to the quantities with indices $12$ in
the above proof. For instance $\phi_{12}$ in (\ref{phiijtest})
equals $\Phi_{k}$ in (\ref{phimaxbis}) as expected from
(\ref{phi}). This allows to check the coherence of our results:
for $n=2$
$$\tr\abs{S_k} = 2 \abs{S^k_{12}} = 2R\frac{\abs{\sin k\pi\omega +\frac{\varphi}{2}}}{\abs{\sin
\omega\pi}}$$ so that 
\begin{equation}
\label{resultatfibre}
d(\xox,\xoz) = 2\pi R\frac{\abs{\sin k\pi\omega
+\frac{\varphi}{2}}}{\abs{\sin \omega\pi}}
\end{equation}
 as expected from
proposition \ref{propcarat}. Indeed when $\tau_0=0$ the domains $\mathcal{T}_\pm$
reduce to the origin and (\ref{defH}) reduces to (\ref{resultatfibre}).

\section{Conclusion and outlook}

Let us conclude on the three generalizations  of  [\citelow{cc}]  presented in this paper: 

- the connected component of the spectral distance for arbitrary dimension $n$ and arbitrary connection
maintains the interesting features of the $n=2$ case. With respect to the horizontal distance, the spectral distance seems to have some transverse properties that could be relevant in  a foliation framework.

- the explicit formula of the distance on all the connected component when $n=2$ shows a non-trivial dependence on $z_\xi$. This is quite unexpected. 

- the explicit formula on the fiber for arbitrary $n$ is remarkably  simple. We deserve its interpretation to further work. Especially it would be interesting
to see whether or not the trace of $\abs{S}$ corresponds to some riemannian distance {\it within} a torus,
as in the $n=2$ case.  In any case the link between the circle and the disk from the spectral metric point of view has to be studied further. Remember that the spectral distance associated to $M_2(\cc)$ only is the euclidean distance on the circle\cite{finite}, while here we find the euclidean distance on the disk. Why considering the circle as a base rather than a point allows the distance to see through the disk ? 

On the circle the difference between the horizontal and the spectral distances is due to the topology of the base. It would be interesting to study local properties in order to see how the curvature enters the game.  The question for sub-riemannian geometry (how to minimize the number of self-intersecting points on a minimal horizontal curve ?) may prevent to go very far in this direction.  Nevertheless one can hope to obtain interesting non trivial results  for base  more complicated than $S^1$ but for which the horizontal distance is still tractable, for instance tori.

Finally, from a more physicist point of view, simple examples of the spectral distance for finite dimensional algebra have been quantized\cite{simple,besnard}.  A similar procedure for covariant Dirac operator would be of great interest.
\newpage

{\bf Acknowledgments:}  Work
partially supported by EU Marie Curie fellowship EIF-025947-QGNC.

\end{document}